\documentclass[11pt]{amsart}
\usepackage{amsmath,amssymb}
\newtheorem{theorem}{Theorem}[section]
\newtheorem{proposition}[theorem]{Proposition}
\newtheorem{corollary}[theorem]{Corollary}
\newtheorem{lemma}[theorem]{Lemma}
\theoremstyle{definition}
\newtheorem{definition}[theorem]{Definition}
\newtheorem{remark}[theorem]{Remark}
\newtheorem{assumption}[theorem]{Assumption}
\newtheorem*{AngleHypothesis}{Matching Angle Hypothesis}

\begin{document}

\title{Scalar curvature rigidity of convex polytopes}
\author{Simon Brendle}
\address{Columbia University \\ 2990 Broadway \\ New York NY 10027 \\ USA}
\begin{abstract}
We prove a scalar curvature rigidity theorem for convex polytopes. The proof uses the Fredholm theory for Dirac operators on manifolds with boundary. A variant of a theorem of Fefferman and Phong plays a central role in our analysis.
\end{abstract}
\thanks{The author is grateful to Yipeng Wang and to the referee for helpful comments on an earlier version of this paper. The author was supported by the National Science Foundation under grant DMS-2103573 and by the Simons Foundation. He acknowledges the hospitality of T\"ubingen University, where part of this work was carried out.}
\maketitle

\section{Introduction}

Let $n \geq 3$ be an integer, and let $\Omega$ be a compact, convex polytope in $\mathbb{R}^n$ with non-empty interior. We may write $\Omega = \bigcap_{i \in I} \{u_i \leq 0\}$, where $u_i$, $i \in I$, is a finite collection of non-constant linear functions defined on $\mathbb{R}^n$. For each $i \in I$, we denote by $N_i \in S^{n-1}$ the outward-pointing unit normal vector to the halfspace $\{u_i \leq 0\}$ with respect to the Euclidean metric.

Let $g$ be a Riemannian metric which is defined on an open set containing $\Omega$. For each $i \in I$, we denote by $\nu_i$ the outward-pointing unit normal vector to the halfspace $\{u_i \leq 0\}$ with respect to the metric $g$. We will assume the following:

\begin{AngleHypothesis}
If $x$ is point in $\partial \Omega$ and $i_1,i_2 \in I$ satisfy $u_{i_1}(x)=u_{i_2}(x)=0$, then $\langle \nu_{i_1},\nu_{i_2} \rangle = \langle N_{i_1},N_{i_2} \rangle$ at the point $x$. Here, the inner product $\langle \nu_{i_1},\nu_{i_2} \rangle$ is computed with respect to the metric $g$, and the inner product $\langle N_{i_1},N_{i_2} \rangle$ is the standard inner product in $\mathbb{R}^n$.
\end{AngleHypothesis}

\begin{theorem}
\label{main.theorem}
Suppose that $n \geq 3$ is an integer, and $\Omega$ is a compact, convex polytope in $\mathbb{R}^n$ with non-empty interior. Let $g$ be a Riemannian metric which is defined on an open set containing $\Omega$ and has nonnegative scalar curvature at each point in $\Omega$. For each $i \in I$, we assume that the mean curvature of the hypersurface $\{u_i=0\}$ with respect to $g$ is nonnegative at each point in $\Omega \cap \{u_i=0\}$. Moreover, we assume that the Matching Angle Hypothesis is satisfied. Then the Riemann curvature tensor of $g$ vanishes at each point in $\Omega$. Moreover, the second fundamental form of the boundary faces of $\Omega$ with respect to $g$ vanishes.
\end{theorem}

Scalar curvature comparison theorems for polytopes were first studied in seminal work of Gromov \cite{Gromov1},\cite{Gromov2},\cite{Gromov3}. In particular, Gromov addressed the case when the dihedral angles are at most $\frac{\pi}{2}$ (see \cite{Gromov2}, Section 3.18). Li \cite{Li1} has used minimal surface techniques to prove a scalar curvature comparison theorem for prisms in dimension $3$. In \cite{Li2}, Li generalized this approach up to dimension $7$. Wang, Xie, and Yu \cite{Wang-Xie-Yu} have proposed a different approach to this problem which is based on the study of Dirac operators on manifolds with corners. 

In this paper, we describe another approach to this problem. As in \cite{Wang-Xie-Yu}, we employ a spinor approach. In contrast to \cite{Wang-Xie-Yu}, we work with boundary value problems for Dirac operators on smooth domains. These types of boundary value problems are well understood thanks to the work of H\"ormander \cite{Hormander2} and B\"ar and Ballmann \cite{Baer-Ballmann1},\cite{Baer-Ballmann2}.

Spinor techniques have long been used in the study of scalar curvature, see e.g. \cite{Bartnik}, \cite{Gromov-Lawson1}, \cite{Gromov-Lawson2}, \cite{Lichnerowicz}, \cite{Llarull}, \cite{Lott}, \cite{Witten}. We refer to the text by Lawson and Michelsohn \cite{Lawson-Michelsohn} for an excellent introduction to spinors and their applications in geometry. 

In the following, we outline the main steps involved in the proof of Theorem \ref{main.theorem}. We approximate a given convex polytope $\Omega$ by a one-parameter family of smooth convex domains $\Omega_\lambda$, where $\lambda$ is assumed to be sufficiently large. The domains $\Omega_\lambda$ form an increasing family of sets, and their union equals the interior of $\Omega$. We consider a sequence $\lambda_l \to \infty$. Let $m = 2^{[\frac{n}{2}]}$ denote the dimension of the space of spinors on flat $\mathbb{R}^n$. For each $l$, we construct a non-trivial $m$-tuple of harmonic spinors $s^{(l)} = (s_1^{(l)},\hdots,s_m^{(l)})$ on the domain $\Omega_{\lambda_l}$ which satisfies a suitable local boundary condition of Lopatinsky-Shapiro-type. To prove the existence of an $m$-tuple of spinors with these properties, we use the Fredholm theory from \cite{Hormander2} together with the deformation invariance of the Fredholm index. We normalize $s^{(l)}$ so that $\int_U \sum_{\alpha=1}^m |s_\alpha^{(l)}|^2 \, d\text{\rm vol}_g = m \, \text{\rm vol}_g(U)$, where $U$ is some fixed Euclidean ball with the property that the closure of $U$ is contained in the interior of $\Omega$. We then apply the Weitzenb\"ock formula to $s^{(l)}$, and integrate over the domain $\Omega_{\lambda_l}$. The resulting integral formula contains a term involving the scalar curvature, as well as a boundary term. Unfortunately, it is not clear if the boundary term has a favorable sign. We are able to control the boundary integral by adapting a deep theorem due to Fefferman and Phong \cite{Fefferman-Phong}. As a result, we are able to show that $\int_{\Omega_{\lambda_l}} \sum_{\alpha=1}^m |\nabla s_\alpha^{(l)}|^2 \, d\text{\rm vol}_g \to 0$ as $l \to \infty$ (see Proposition \ref{covariant.derivative.of.s.small.in.L2} below). By passing to the limit as $l \to \infty$, we obtain an $m$-tuple of parallel spinors $s = (s_1,\hdots,s_m)$ which is defined in the interior of $\Omega$. In particular, $\langle s_\alpha,s_\beta \rangle = z_{\alpha\beta}$ for some fixed matrix $z \in \text{\rm End}(\mathbb{C}^m)$. Finally, by exploiting the boundary condition, we show that $z$ is the identity. In other words, $s_1,\hdots,s_m$ are orthonormal at each point in the interior of $\Omega$. As a consequence, the Riemann curvature tensor of $g$ vanishes identically.

\section{A boundary value problem for the Dirac operator on an odd-dimensional domain with smooth boundary}

Throughout this section, we assume that $n \geq 3$ is an odd integer. Let $\{E_1,\hdots,E_n\}$ denote the standard basis of $\mathbb{R}^n$. Let $\text{\rm Cl}(n,\mathbb{C})$ denote the Clifford algebra. The spin representation gives a surjective algebra homomorphism $\hat{\rho}: \text{\rm Cl}(n,\mathbb{C}) \to \text{\rm End}(\Delta_n)$, where $\Delta_n$ is a complex vector space of dimension $m = 2^{[\frac{n}{2}]}$ equipped with a Hermitian inner product. For each $a=1,\hdots,n$, the map $\hat{\rho}(E_a) \in \text{\rm End}(\Delta_n)$ is skew-adjoint. Moreover, the Clifford relations 
\begin{equation} 
\label{clifford.relations}
\hat{\rho}(E_a) \hat{\rho}(E_b) + \hat{\rho}(E_b) \hat{\rho}(E_a) = -2\delta_{ab} \, \text{\rm id} 
\end{equation}
hold for all $a,b=1,\hdots,n$. Since $n$ is odd, the product $E_1 \cdots E_n \in \text{\rm Cl}(n,\mathbb{C})$ commutes with every element of $\text{\rm Cl}(n,\mathbb{C})$. Since $\hat{\rho}$ is surjective, it follows that $\hat{\rho}(E_1) \cdots \hat{\rho}(E_n) \in \text{\rm End}(\Delta_n)$ commutes with every element of $\text{\rm End}(\Delta_n)$. Therefore, $\hat{\rho}(E_1) \cdots \hat{\rho}(E_n) \in \text{\rm End}(\Delta_n)$ is a scalar multiple of the identity. It is straightforward to see that 
\begin{equation} 
\label{volume.form}
i^{\frac{n+1}{2}} \, \hat{\rho}(E_1) \cdots \hat{\rho}(E_n) = \pm \text{\rm id}. 
\end{equation} 
The sign in (\ref{volume.form}) depends on the choice of $\hat{\rho}$ (see \cite{Lawson-Michelsohn}, Proposition 5.9). In the following, we assume that $\hat{\rho}$ is chosen so that $i^{\frac{n+1}{2}} \, \hat{\rho}(E_1) \cdots \hat{\rho}(E_n) = \text{\rm id}$. 

Let us fix an orthonormal basis $\{\hat{s}_1,\hdots,\hat{s}_m\}$ of $\Delta_n$. We define 
\begin{equation} 
\omega_{a\alpha\beta} = \langle \hat{\rho}(E_a) \, \hat{s}_\alpha,\hat{s}_\beta \rangle 
\end{equation} 
for $a=1,\hdots,n$ and $\alpha,\beta=1,\hdots,m$. The matrices $\omega_1,\hdots,\omega_n \in \text{\rm End}(\mathbb{C}^m)$ are skew-Hermitian and satisfy $\omega_a \omega_b + \omega_b \omega_a = -2\delta_{ab} \, \text{\rm id}$ for all $a,b=1,\hdots,n$. Moreover, $i^{\frac{n+1}{2}} \, \omega_1 \cdots \omega_n = \text{\rm id}$ by our choice of $\hat{\rho}$. Finally, 
\begin{equation} 
\label{surjectivity.of.hat.rho}
\text{\rm End}(\mathbb{C}^m) = \text{\rm span}\{\omega_{a_1} \cdots \omega_{a_r}: 1 \leq a_1 < \hdots < a_r \leq n\} 
\end{equation}
since $\hat{\rho}$ is surjective.

\begin{lemma}
\label{algebraic.fact}
Assume that $n \geq 3$ is an odd integer. If $z \in \text{\rm End}(\mathbb{C}^m)$ anti-commutes with $\omega_a \in \text{\rm End}(\mathbb{C}^m)$ for each $a=1,\hdots,n$, then $z=0$. 
\end{lemma}

\textbf{Proof.} 
Suppose that $z \in \text{\rm End}(\mathbb{C}^m)$ anti-commutes with $\omega_a \in \text{\rm End}(\mathbb{C}^m)$ for each $a=1,\hdots,n$. Since $n$ is odd, it follows that $z$ anti-commutes with the product $\omega_1 \cdots \omega_n \in \text{\rm End}(\mathbb{C}^m)$. Since $\omega_1 \cdots \omega_n$ is a non-zero multiple of the identity, we conclude that $z=0$. This completes the proof of Lemma \ref{algebraic.fact}. \\

In the remainder of this section, we assume that $\Omega$ is a compact domain in $\mathbb{R}^n$ with smooth boundary $\partial \Omega = \Sigma$. Let $g$ be a Riemannian metric which is defined on an open set containing $\Omega$. We denote by $\nu$ the outward-pointing unit normal vector field with respect to the metric $g$. We denote by $H$ the mean curvature of $\Sigma$ with respect to $g$, defined as the trace of the fundamental form of $\Sigma$. Under our sign convention for $H$, the mean curvature vector of $\Sigma$ is given by $-H\nu$.

Let $\mathcal{S}$ denote the spinor bundle over $\Omega$ with respect to the metric $g$. Note that $\mathcal{S}$ is a complex vector bundle of rank $m$ equipped with a Hermitian inner product. As a bundle, we may identify $\mathcal{S}$ with the trivial bundle $\Omega \times \Delta_n \to \Omega$. Each tangent vector $\xi \in T_x \Omega$ induces a skew-adjoint map $\rho(\xi) \in \text{\rm End}(\mathcal{S}_x)$. If $\{e_1,\hdots,e_n\}$ is an orthonormal frame with respect to the metric $g$, then 
\[\rho(e_k) \rho(e_l) + \rho(e_k) \rho(e_l) = -2\delta_{kl} \, \text{\rm id}\] 
for $k,l=1,\hdots,n$. Moreover, if $\{e_1,\hdots,e_n\}$ is a positively oriented orthonormal frame with respect to the metric $g$, then $i^{\frac{n+1}{2}} \, \rho(e_1) \cdots \rho(e_n) = \text{\rm id}$. For abbreviation, we write $\xi \cdot s$ instead of $\rho(\xi) \, s$. This is referred to as Clifford multiplication.

We next consider the spin connection with respect to the metric $g$. The spin connection is a connection $\nabla$ on $\mathcal{S}$ which is compatible with Clifford multiplication and which is compatible with the Hermitian inner product on $\mathcal{S}$. The Dirac operator is defined by 
\[\mathcal{D} s = \sum_{k=1}^n e_k \cdot \nabla_{e_k} s,\] 
where $s$ denotes a section of $\mathcal{S}$ and $\{e_1,\hdots,e_n\}$ is a local orthonormal frame on $\Omega$. The boundary Dirac operator $\mathcal{D}^\Sigma$ is given by  
\[\mathcal{D}^\Sigma s = \sum_{k=1}^{n-1} \nu \cdot e_k \cdot \nabla_{e_k} s + \frac{1}{2} \, H \, s,\] 
where $s$ denotes a section of $\mathcal{S}|_\Sigma$ and $\{e_1,\hdots,e_{n-1}\}$ is a local orthonormal frame on $\Sigma$. Note that $\mathcal{D}^\Sigma$ is formally self-adjoint.

In the following, we will consider the Dirac operator acting on $m$-tuples of spinors. To fix notation, we define a complex vector bundle $\mathcal{E}$ over $\Omega$ by 
\[\mathcal{E} = \underbrace{\mathcal{S} \oplus \hdots \oplus \mathcal{S}}_{\text{\rm $m$ times}}.\] 
Note that $\mathcal{E}$ has rank $m^2$. A section of $\mathcal{E}$ can be identified with an $m$-tuple of spinors $s = (s_1,\hdots,s_m)$ defined on $\Omega$. We may view $\mathcal{D}$ as an operator acting on sections of $\mathcal{E}$, and we may view $\mathcal{D}^\Sigma$ as an operator acting on sections of $\mathcal{E}|_\Sigma$.

In the next step, we introduce a local boundary condition of Lopatinsky-Shapiro-type. To formulate the boundary condition, we assume that a smooth map $N: \Sigma \to S^{n-1}$ is given. 

\begin{definition}
\label{definition.of.chi}
We define a bundle map $\chi: \mathcal{E}|_\Sigma \to \mathcal{E}|_\Sigma$ by 
\[(\chi s)_\alpha = -\sum_{a=1}^n \sum_{\beta=1}^m \langle N,E_a \rangle \, \omega_{a\alpha\beta} \, \nu \cdot s_\beta\] 
for $\alpha=1,\hdots,m$. Moreover, we define a bundle map $\mathcal{B}: \mathcal{E}|_\Sigma \to \mathcal{E}|_\Sigma$ by 
\[(\mathcal{B} s)_\alpha = \sum_{k=1}^{n-1} \sum_{a=1}^n \sum_{\beta=1}^m \langle dN(e_k),E_a \rangle \, \omega_{a\alpha\beta} \, e_k \cdot s_\beta\] 
for $\alpha=1,\hdots,m$. Here, $\{e_1,\hdots,e_{n-1}\}$ is a local orthonormal frame on $\Sigma$.
\end{definition}

\begin{lemma}
\label{properties.of.chi}
The map $\chi$ is self-adjoint. Moreover, $\chi^2 = \text{\rm id}$. 
\end{lemma} 

\textbf{Proof.} 
Suppose that $s = (s_1,\hdots,s_m)$ and $t = (t_1,\hdots,t_m)$ are two $m$-tuples of spinors. We compute 
\begin{align*} 
\sum_{\alpha=1}^m \langle (\chi s)_\alpha,t_\alpha \rangle 
&= -\sum_{a=1}^n \sum_{\alpha,\beta=1}^m \langle N,E_a \rangle \, \omega_{a\alpha\beta} \, \langle \nu \cdot s_\beta,t_\alpha \rangle \\ 
&= -\sum_{a=1}^n \sum_{\alpha,\beta=1}^m \langle N,E_a \rangle \, \overline{\omega_{a\beta\alpha}} \, \langle s_\beta,\nu \cdot t_\alpha \rangle \\ 
&= \sum_{\beta=1}^m \langle s_\beta,(\chi t)_\beta \rangle. 
\end{align*}
Moreover, 
\[(\chi^2 s)_\alpha = -\sum_{a,b=1}^n \sum_{\beta,\gamma=1}^m \langle N,E_a \rangle \, \langle N,E_b \rangle \, \omega_{a\alpha\beta} \, \omega_{b\beta\gamma} \, s_\gamma = s_\alpha\] 
for $\alpha=1,\hdots,m$. This completes the proof of Lemma \ref{properties.of.chi}. \\

\begin{lemma} 
\label{properties.of.B}
The map $\mathcal{B}$ is self-adjoint. Moreover, $\chi$ and $\mathcal{B}$ commute. 
\end{lemma}

\textbf{Proof.} 
Suppose that $s = (s_1,\hdots,s_m)$ and $t = (t_1,\hdots,t_m)$ are two $m$-tuples of spinors. Let $\{e_1,\hdots,e_{n-1}\}$ be a local orthonormal frame on $\Sigma$. Then 
\begin{align*} 
\sum_{\alpha=1}^m \langle (\mathcal{B} s)_\alpha,t_\alpha \rangle 
&= \sum_{k=1}^{n-1} \sum_{a=1}^n \sum_{\alpha,\beta=1}^m \langle dN(e_k),E_a \rangle \, \omega_{a\alpha\beta} \, \langle e_k \cdot s_\beta,t_\alpha \rangle \\ 
&= \sum_{k=1}^{n-1} \sum_{a=1}^n \sum_{\alpha,\beta=1}^m \langle dN(e_k),E_a \rangle \, \overline{\omega_{a\beta\alpha}} \, \langle s_\beta,e_k \cdot t_\alpha \rangle \\ 
&= \sum_{\beta=1}^m \langle s_\beta,(\mathcal{B} t)_\beta \rangle. 
\end{align*} 
This shows that $\mathcal{B}$ is self-adjoint. Moreover, 
\begin{align*} 
&(\chi \mathcal{B} s)_\alpha - (\mathcal{B} \chi s)_\alpha \\ 
&= -\sum_{k=1}^{n-1} \sum_{a,b=1}^n \sum_{\beta,\gamma=1}^m \langle N,E_a \rangle \, \langle dN(e_k),E_b \rangle \, \omega_{a\alpha\beta} \, \omega_{b\beta\gamma} \, \nu \cdot e_k \cdot s_\gamma \\ 
&+ \sum_{k=1}^{n-1} \sum_{a,b=1}^n \sum_{\beta,\gamma=1}^m \langle dN(e_k),E_a \rangle \, \langle N,E_b \rangle \, \omega_{a\alpha\beta} \, \omega_{b\beta\gamma} \, e_k \cdot \nu \cdot s_\gamma \\ 
&= -\sum_{k=1}^{n-1} \sum_{a,b=1}^n \sum_{\beta,\gamma=1}^m \langle N,E_a \rangle \, \langle dN(e_k),E_b \rangle \, (\omega_{a\alpha\beta} \, \omega_{b\beta\gamma} + \omega_{b\alpha\beta} \, \omega_{a\beta\gamma}) \, \nu \cdot e_k \cdot s_\gamma \\ 
&= 2 \sum_{k=1}^{n-1} \langle N,dN(e_k) \rangle \, \nu \cdot e_k \cdot s_\alpha \\ 
&= 0 
\end{align*} 
for $\alpha=1,\hdots,m$. Thus, $\chi$ and $\mathcal{B}$ commute. This completes the proof of Lemma \ref{properties.of.B}. \\

\begin{proposition}
\label{commutator.of.chi.and.boundary.Dirac.operator}
Suppose that $s = (s_1,\hdots,s_m)$ is an $m$-tuple of spinors. Then 
\[\chi \mathcal{D}^\Sigma s + \mathcal{D}^\Sigma \chi s = -\mathcal{B} s.\] 
\end{proposition} 

\textbf{Proof.} 
Let $\{e_1,\hdots,e_{n-1}\}$ denote a local orthonormal frame on $\Sigma$. We compute 
\begin{align*} 
(\chi \mathcal{D}^\Sigma s)_\alpha 
&= -\sum_{a=1}^n \sum_{\beta=1}^m \langle N,E_a \rangle \, \omega_{a\alpha\beta} \, \nu \cdot (\mathcal{D}^\Sigma s_\beta) \\ 
&= \sum_{k=1}^{n-1} \sum_{a=1}^n \sum_{\beta=1}^m \langle N,E_a \rangle \, \omega_{a\alpha\beta} \, e_k \cdot \nabla_{e_k} s_\beta + \frac{1}{2} \, H \, (\chi s)_\alpha 
\end{align*} 
and 
\begin{align*} 
(\mathcal{D}^\Sigma \chi s)_\alpha 
&= \sum_{k=1}^{n-1} \nu \cdot e_k \cdot \nabla_{e_k} (\chi s)_\alpha + \frac{1}{2} \, H \, (\chi s)_\alpha \\ 
&= -\sum_{k=1}^{n-1} e_k \cdot \nabla_{e_k} (\nu \cdot \chi s)_\alpha - \frac{1}{2} \, H \, (\chi s)_\alpha \\ 
&= -\sum_{k=1}^{n-1} \sum_{a=1}^n \sum_{\beta=1}^m \langle dN(e_k),E_a \rangle \, \omega_{a\alpha\beta} \, e_k \cdot s_\beta \\ 
&-\sum_{k=1}^{n-1} \sum_{a=1}^n \sum_{\beta=1}^m \langle N,E_a \rangle \, \omega_{a\alpha\beta} \, e_k \cdot \nabla_{e_k} s_\beta - \frac{1}{2} \, H \, (\chi s)_\alpha 
\end{align*}
for $\alpha=1,\hdots,m$. Putting these facts together, the assertion follows. This completes the proof of Proposition \ref{commutator.of.chi.and.boundary.Dirac.operator}. \\

\begin{corollary}
\label{commutator.of.chi.and.A}
Let $\mathcal{A} = \mathcal{D}^\Sigma + \frac{1}{2} \chi \mathcal{B}$. Then $\mathcal{A}$ is formally self-adjoint. Moreover, $\mathcal{A}$ anti-commutes with $\chi$. 
\end{corollary}

\textbf{Proof.} 
It follows from Lemma \ref{properties.of.chi} and Lemma \ref{properties.of.B} that $\chi \mathcal{B}$ is self-adjoint. Since $\mathcal{D}^\Sigma$ is formally self-adjoint, we conclude that $\mathcal{A}$ is formally self-adjoint. This proves the first statement. To prove the second statement, suppose that $s = (s_1,\hdots,s_m)$ is an $m$-tuple of spinors. Using Lemma \ref{properties.of.chi}, Lemma \ref{properties.of.B}, and Proposition \ref{commutator.of.chi.and.boundary.Dirac.operator}, we obtain 
\begin{align*} 
\chi \mathcal{A} s + \mathcal{A} \chi s 
&= \chi \mathcal{D}^\Sigma s + \mathcal{D}^\Sigma \chi s + \frac{1}{2} \chi \chi \mathcal{B} s + \frac{1}{2} \chi \mathcal{B} \chi s \\ 
&=\chi \mathcal{D}^\Sigma s +  \mathcal{D}^\Sigma \chi s + \mathcal{B} s \\ 
&= 0. 
\end{align*}
This completes the proof of Corollary \ref{commutator.of.chi.and.A}. \\

At this point, we recall a definition from linear algebra (see e.g. \cite{Zhan}, p.~92). 

\begin{definition}
Let $V$ and $W$ be finite-dimensional real vector spaces of the same dimension, each of them equipped with an inner product. Let $L: V \to W$ be a linear map. The trace norm of $L$ is defined by $\|L\|_{\text{\rm tr}} = \sup_Q \text{\rm tr}(QL)$, where the supremum is taken over all linear isometries $Q: W \to V$. Equivalently, $\|L\|_{\text{\rm tr}}$ can be characterized as the sum of the singular values of $L$.
\end{definition}

It is easy to see from the definition that the trace norm satisfies the triangle inequality. \\

\begin{lemma} 
\label{bound.for.B}
Suppose that $s = (s_1,\hdots,s_m)$ is an $m$-tuple of spinors. Then 
\[\bigg | \sum_{\alpha=1}^m \langle (\mathcal{B} s)_\alpha,s_\alpha \rangle \bigg | \leq \|dN\|_{\text{\rm tr}} \, \bigg ( \sum_{\alpha=1}^m |s_\alpha|^2 \bigg )\] 
at each point $x \in \Sigma$. Here, $\|dN\|_{\text{\rm tr}}$ denotes the trace norm of the differential $dN: T_x \Sigma \to T_{N(x)} S^{n-1}$. The tangent space $T_x \Sigma$ is equipped with the restriction of the inner product $g$, and the tangent space $T_{N(x)} S^{n-1}$ is equipped with the restriction of the standard inner product on $\mathbb{R}^n$.
\end{lemma} 

\textbf{Proof.} 
Fix a point $x \in \Sigma$. Let $\lambda_1,\hdots,\lambda_{n-1} \geq 0$ denote the singular values of the differential $dN: T_x \Sigma \to T_{N(x)} S^{n-1}$. We can find an orthonormal basis $\{e_1,\hdots,e_{n-1}\}$ of $T_x \Sigma$ and an orthonormal basis $\{\hat{E}_1,\hdots,\hat{E}_{n-1}\}$ of $T_{N(x)} S^{n-1}$ such that $dN(e_k) = \lambda_k \, \hat{E}_k$ for each $k=1,\hdots,n-1$. Then 
\begin{align*} 
&\sum_{\alpha=1}^m \bigg | \sum_{a=1}^n \sum_{\beta=1}^m \langle \hat{E}_k,E_a \rangle \, \omega_{a\alpha\beta} \, e_k \cdot s_\beta \bigg |^2 \\ 
&= -\sum_{a,b=1}^n \sum_{\alpha,\beta,\gamma=1}^m \langle \hat{E}_k,E_a \rangle \, \langle \hat{E}_k,E_b \rangle \, \omega_{a\alpha\beta} \, \omega_{b\gamma\alpha} \, \langle s_\beta,s_\gamma \rangle \\ 
&= -\frac{1}{2} \sum_{a,b=1}^n \sum_{\alpha,\beta,\gamma=1}^m \langle \hat{E}_k,E_a \rangle \, \langle \hat{E}_k,E_b \rangle \, (\omega_{a\gamma\alpha} \, \omega_{b\alpha\beta} + \omega_{b\gamma\alpha} \,  \omega_{a\alpha\beta}) \, \langle s_\beta,s_\gamma \rangle \\ 
&= \sum_{\alpha=1}^m |s_\alpha|^2 
\end{align*} 
for each $k=1,\hdots,n-1$. Using the Cauchy-Schwarz inequality, we obtain 
\begin{align*} 
&\bigg | \sum_{a=1}^n \sum_{\alpha,\beta=1}^m \langle \hat{E}_k,E_a \rangle \, \omega_{a\alpha\beta} \, \langle e_k \cdot s_\beta,s_\alpha \rangle \bigg | \\ 
&\leq \bigg ( \sum_{\alpha=1}^m \bigg | \sum_{a=1}^n \sum_{\beta=1}^m \langle \hat{E}_k,E_a \rangle \, \omega_{a\alpha\beta} \, e_k \cdot s_\beta \bigg |^2 \bigg )^{\frac{1}{2}} \, \bigg ( \sum_{\alpha=1}^m |s_\alpha|^2 \bigg )^{\frac{1}{2}} = \sum_{\alpha=1}^m |s_\alpha|^2 
\end{align*}
for each $k=1,\hdots,n-1$. Summation over $k=1,\hdots,n-1$ gives 
\begin{align*} 
\bigg | \sum_{\alpha=1}^m \langle (\mathcal{B} s)_\alpha,s_\alpha \rangle \bigg | 
&= \bigg | \sum_{k=1}^{n-1} \sum_{a=1}^n \sum_{\alpha,\beta=1}^m \langle dN(e_k),E_a \rangle \, \omega_{a\alpha\beta} \, \langle e_k \cdot s_\beta,s_\alpha \rangle \bigg | \\ 
&= \bigg | \sum_{k=1}^{n-1} \lambda_k \, \bigg ( \sum_{a=1}^n \sum_{\alpha,\beta=1}^m \langle \hat{E}_k,E_a \rangle \, \omega_{a\alpha\beta} \, \langle e_k \cdot s_\beta,s_\alpha \rangle \bigg ) \bigg | \\ 
&\leq \bigg ( \sum_{k=1}^{n-1} \lambda_k \bigg ) \, \bigg ( \sum_{\alpha=1}^m |s_\alpha|^2 \bigg ). 
\end{align*}  
This completes the proof of Lemma \ref{bound.for.B}. \\

\begin{proposition}
\label{main.inequality} 
Suppose that $s = (s_1,\hdots,s_m)$ is an $m$-tuple of spinors. Then 
\begin{align*} 
&-\int_\Omega \sum_{\alpha=1}^m |\mathcal{D} s_\alpha|^2 \, d\text{\rm vol}_g + \int_\Omega \sum_{\alpha=1}^m |\nabla s_\alpha|^2 \, d\text{\rm vol}_g + \frac{1}{4} \int_\Omega \sum_{\alpha=1}^m R \, |s_\alpha|^2 \, d\text{\rm vol}_g \\ 
&\leq \frac{1}{2} \int_\Sigma \sum_{\alpha=1}^m \langle \mathcal{D}^\Sigma s_\alpha,s_\alpha - (\chi s)_\alpha \rangle \, d\sigma_g + \frac{1}{2} \int_\Sigma \sum_{\alpha=1}^m \langle s_\alpha - (\chi s)_\alpha,\mathcal{D}^\Sigma s_\alpha \rangle \, d\sigma_g \\ 
&- \frac{1}{2} \int_\Sigma (H - \|dN\|_{\text{\rm tr}}) \, \bigg ( \sum_{\alpha=1}^m |s_\alpha|^2 \bigg ) \, d\sigma_g. 
\end{align*} 
\end{proposition}

\textbf{Proof.}
By the Weitzenb\"ock formula, $\mathcal{D}^2 s_\alpha = -\Delta s_\alpha + \frac{1}{4} \, R \, s_\alpha$, where $\Delta$ denotes the connection Laplacian on the spinor bundle. Using the divergence theorem, we obtain 
\begin{align*} 
&-\int_\Omega \sum_{\alpha=1}^m |\mathcal{D} s_\alpha|^2 \, d\text{\rm vol}_g + \int_\Omega \sum_{\alpha=1}^m |\nabla s_\alpha|^2 \, d\text{\rm vol}_g + \frac{1}{4} \int_\Omega \sum_{\alpha=1}^m R \, |s_\alpha|^2 \, d\text{\rm vol}_g \\ 
&= \int_\Sigma \sum_{\alpha=1}^m \langle \nu \cdot \mathcal{D} s_\alpha,s_\alpha \rangle \, d\sigma_g + \int_\Sigma \sum_{\alpha=1}^m \langle \nabla_\nu s_\alpha,s_\alpha \rangle \, d\sigma_g. 
\end{align*} 
Note that $\nu \cdot \mathcal{D} s_\alpha + \nabla_\nu s_\alpha = \mathcal{D}^\Sigma s_\alpha - \frac{1}{2} \, H \, s_\alpha$ at each point on $\Sigma$. This implies 
\begin{align*} 
&-\int_\Omega \sum_{\alpha=1}^m |\mathcal{D} s_\alpha|^2 \, d\text{\rm vol}_g + \int_\Omega \sum_{\alpha=1}^m |\nabla s_\alpha|^2 \, d\text{\rm vol}_g + \frac{1}{4} \int_\Omega \sum_{\alpha=1}^m R \, |s_\alpha|^2 \, d\text{\rm vol}_g \\ 
&= \int_\Sigma \langle \mathcal{D}^\Sigma s_\alpha,s_\alpha \rangle \, d\sigma_g - \frac{1}{2} \int_\Sigma \sum_{\alpha=1}^m H \, |s_\alpha|^2 \, d\sigma_g.
\end{align*} 
On the other hand, using the fact that $\chi$ is self-adjoint and $\mathcal{D}^\Sigma$ is formally self-adjoint, we obtain 
\begin{align*} 
&\int_\Sigma \sum_{\alpha=1}^m \langle \mathcal{D}^\Sigma s_\alpha,(\chi s)_\alpha \rangle \, d\sigma_g + \int_\Sigma \sum_{\alpha=1}^m \langle (\chi s)_\alpha,\mathcal{D}^\Sigma s_\alpha \rangle \, d\sigma_g \\
&= \int_\Sigma \sum_{\alpha=1}^m \langle (\chi \mathcal{D}^\Sigma s)_\alpha,s_\alpha \rangle \, d\sigma_g + \int_\Sigma \sum_{\alpha=1}^m \langle (\mathcal{D}^\Sigma \chi s)_\alpha,s_\alpha \rangle \, d\sigma_g \\ 
&= -\int_\Sigma \sum_{\alpha=1}^m \langle (\mathcal{B} s)_\alpha,s_\alpha \rangle \, d\sigma_g.
\end{align*}
In the last step, we have used Proposition \ref{commutator.of.chi.and.boundary.Dirac.operator}. Putting these facts together, we conclude that 
\begin{align*} 
&-\int_\Omega \sum_{\alpha=1}^m |\mathcal{D} s_\alpha|^2 \, d\text{\rm vol}_g + \int_\Omega \sum_{\alpha=1}^m |\nabla s_\alpha|^2 \, d\text{\rm vol}_g + \frac{1}{4} \int_\Omega \sum_{\alpha=1}^m R \, |s_\alpha|^2 \, d\text{\rm vol}_g \\ 
&= \frac{1}{2} \int_\Sigma \sum_{\alpha=1}^m \langle \mathcal{D}^\Sigma s_\alpha,s_\alpha \rangle \, d\sigma_g + \frac{1}{2} \int_\Sigma \sum_{\alpha=1}^m \langle s_\alpha,\mathcal{D}^\Sigma s_\alpha \rangle \, d\sigma_g - \frac{1}{2} \int_\Sigma \sum_{\alpha=1}^m H \, |s_\alpha|^2 \, d\sigma_g \\ 
&= \frac{1}{2} \int_\Sigma \sum_{\alpha=1}^m \langle \mathcal{D}^\Sigma s_\alpha,s_\alpha - (\chi s)_\alpha \rangle \, d\sigma_g + \frac{1}{2} \int_\Sigma \sum_{\alpha=1}^m \langle s_\alpha - (\chi s)_\alpha,\mathcal{D}^\Sigma s_\alpha \rangle \, d\sigma_g \\ 
&- \frac{1}{2} \int_\Sigma \sum_{\alpha=1}^m \langle (\mathcal{B} s)_\alpha,s_\alpha \rangle \, d\sigma_g - \frac{1}{2} \int_\Sigma \sum_{\alpha=1}^m H \, |s_\alpha|^2 \, d\sigma_g. 
\end{align*} 
Hence, the assertion follows from Lemma \ref{bound.for.B}. \\

\begin{corollary} 
\label{rigidity.1}
Assume that $R \geq 0$ at each point in $\Omega$ and $H \geq \|dN\|_{\text{\rm tr}}$ at each point on $\Sigma$. Suppose that $s = (s_1,\hdots,s_m)$ is an $m$-tuple of harmonic spinors on $\Omega$ which satisfies the boundary condition $\chi s = s$ at each point on $\Sigma$. Then $s$ is parallel.
\end{corollary}

Replacing $N$ by $-N$, we can draw the following conclusion: 

\begin{corollary}
\label{rigidity.2} 
Assume that $R \geq 0$ at each point in $\Omega$ and $H \geq \|dN\|_{\text{\rm tr}}$ at each point on $\Sigma$. Suppose that $s = (s_1,\hdots,s_m)$ is an $m$-tuple of harmonic spinors on $\Omega$ which satisfies the boundary condition $\chi s = -s$ at each point on $\Sigma$. Then $s$ is parallel.
\end{corollary}

\begin{lemma}
\label{clifford.multiplication.anti-commutes.with.chi}
Assume that $x \in \Sigma$ and $\xi \in T_x \Sigma$. Then the linear map $(s_1,\hdots,s_m) \mapsto (i \, \nu \cdot \xi \cdot s_1,\hdots,i \, \nu \cdot \xi \cdot s_m)$ anti-commutes with $\chi$. In particular, $\dim \ker (\text{\rm id}-\chi) = \dim \ker (\text{\rm id}+\chi)$.
\end{lemma}

\textbf{Proof.} 
This follows immediately from the definition of $\chi$. \\

\begin{definition}
We denote by $\mathcal{F} = \mathcal{E}|_\Sigma$ the restriction of $\mathcal{E}$ to $\Sigma$. Moreover, we write $\mathcal{F} = \mathcal{F}^+ \oplus \mathcal{F}^-$, where $\mathcal{F}^+ = \text{\rm ker}(\text{\rm id}-\chi)$ and $\mathcal{F}^- = \text{\rm ker}(\text{\rm id}+\chi)$. 
\end{definition}

Note that $\mathcal{F}$ is a complex vector bundle over $\Sigma$ of rank $m^2$. It follows from Lemma \ref{clifford.multiplication.anti-commutes.with.chi} that $\mathcal{F}^+$ and $\mathcal{F}^-$ are complex subbundles of $\mathcal{F}$ of rank $\frac{m^2}{2}$.

\begin{proposition}
\label{fredholm}
Suppose that $\Omega$ is a compact domain in $\mathbb{R}^n$ with smooth boundary $\partial \Omega = \Sigma$. Let $g$ be a Riemannian metric which is defined on an open set containing $\Omega$, and let $N: \Sigma \to S^{n-1}$ be a smooth map. Then the operator 
\[H^1(\Omega,\mathcal{E}) \to L^2(\Omega,\mathcal{E}) \oplus H^{\frac{1}{2}}(\Sigma,\mathcal{F}^-), \quad s \mapsto (\mathcal{D} s,s-\chi s)\] 
is a Fredholm operator. The kernel of this operator is a finite-dimensional subspace of $C^\infty(\Omega,\mathcal{E})$. The range of this operator is defined by finitely many $C^\infty$ relations.
\end{proposition}

\textbf{Proof.} 
We will show that the boundary value problem is elliptic in the sense of Definition 20.1.1 in H\"ormander's book \cite{Hormander2}. To that end, we fix a point $x \in \Sigma$. Moreover, we consider a vector $\xi \in T_x \Omega$ with the property that $\xi$ is not a scalar multiple of $\nu$. Following H\"ormander \cite{Hormander2}, we denote by $M_{x,\xi}^+$ the set of all functions $u: \mathbb{R} \to \mathcal{F}_x$ which solve the linear ODE 
\[i \, \xi \cdot u(t) - \nu \cdot \, \frac{d}{dt} u(t) = 0\] 
and which are bounded on the interval $[0,\infty)$. (Note that $D_t = -i \frac{d}{dt}$ in H\"ormander's notation; see \cite{Hormander1}, p.~160.) If we fix a real number $a$, then the function $t \mapsto u(t)$ belongs to the space $M_{x,\xi}^+$ if and only if the function $t \mapsto e^{iat} \, u(t)$ belongs to the space $M_{x,\xi+a\nu}^+$.

We claim that the linear map 
\[M_{x,\xi}^+ \to \mathcal{F}_x^-, \quad u \mapsto u(0)-\chi u(0)\] 
is bijective. If the claim is true for some vector $\xi \in T_x \Omega$, then the claim is also true for $\xi+a\nu$, where $a$ is an arbitrary real number. Hence, it suffices to verify the claim in the special case when $\xi \in T_x \Sigma$ and $\xi \neq 0$. 

In the following, we assume that $\xi \in T_x \Sigma$ and $\xi \neq 0$. We define a linear map $L: \mathcal{F}_x \to \mathcal{F}_x$ by $Ls = i \, \nu \cdot \xi \cdot s$. Note that $L$ is self-adjoint and $L^2 = |\xi|^2 \, \text{\rm id}$. Therefore, 
\[\mathcal{F}_x = \ker (|\xi| \, \text{\rm id}-L) \oplus \ker (|\xi| \, \text{\rm id}+L).\] 
The space $M_{x,\xi}^+$ consists of all functions of the form $t \mapsto e^{-t|\xi|} \, u_0$, where $u_0 \in \ker(|\xi| \, \text{\rm id}-L)$. 

By Lemma \ref{clifford.multiplication.anti-commutes.with.chi}, $L$ anti-commutes with $\chi$. This implies 
\begin{equation} 
\label{eigenspaces.of.L.have.the.same.dimension}
\dim \ker (|\xi| \, \text{\rm id}-L) = \dim \ker(|\xi| \, \text{\rm id}+L) 
\end{equation} 
and 
\begin{equation} 
\label{intersection.of.eigenspace.of.L.with.eigenspace.of.chi}
\ker (|\xi| \, \text{\rm id}-L) \cap \ker(\text{\rm id}-\chi) = \{0\}. 
\end{equation} 
Using (\ref{eigenspaces.of.L.have.the.same.dimension}) and Lemma \ref{clifford.multiplication.anti-commutes.with.chi}, we obtain 
\[\dim M_{x,\xi}^+ = \dim \ker(|\xi| \, \text{\rm id}-L) = \frac{m^2}{2} = \dim \mathcal{F}_x^-.\] 
Moreover, it follows from (\ref{intersection.of.eigenspace.of.L.with.eigenspace.of.chi}) that the linear map 
\[M_{x,\xi}^+ \to \mathcal{F}_x^-, \quad u \mapsto u(0)-\chi u(0)\] 
has trivial kernel. Therefore, the latter map is bijective, as claimed. 

To summarize, we have shown that the boundary value problem is elliptic in the sense of Definition 20.1.1 in \cite{Hormander2}. The assertion now follows from Theorem 20.1.2 and Theorem 20.1.8 in \cite{Hormander2}. Note that the function spaces appearing in (20.1.2) in \cite{Hormander2} are the usual Sobolev spaces; see Definition 7.9.1 in \cite{Hormander1} and Section B.2 in \cite{Hormander2} for the relevant definitions. This completes the proof of Proposition \ref{fredholm}. \\

\begin{proposition}
\label{index}
Assume that $n \geq 3$ is an odd integer. Suppose that $\Omega$ is a compact, convex domain in $\mathbb{R}^n$ with smooth boundary $\partial \Omega = \Sigma$. Let $g$ be a Riemannian metric which is defined on an open set containing $\Omega$, and let $N: \Sigma \to S^{n-1}$ be a smooth map which is homotopic to the Euclidean Gauss map of $\Sigma$. Then the operator 
\[H^1(\Omega,\mathcal{E}) \to L^2(\Omega,\mathcal{E}) \oplus H^{\frac{1}{2}}(\Sigma,\mathcal{F}^-), \quad s \mapsto (\mathcal{D} s,s-\chi s)\] 
has Fredholm index at least $1$. 
\end{proposition}

\textbf{Proof.} 
We first consider the special case when $g$ is the Euclidean metric and $N$ is the Euclidean Gauss map of $\Sigma$. In this case, a spinor can be viewed as a smooth function taking values in $\Delta_n$. We claim that the kernel of the operator 
\[H^1(\Omega,\mathcal{E}) \to L^2(\Omega,\mathcal{E}) \oplus H^{\frac{1}{2}}(\Sigma,\mathcal{F}^-), \quad s \mapsto (\mathcal{D} s,s-\chi s)\] 
has dimension at least $1$. To see this, recall that $\{\hat{s}_1,\hdots,\hat{s}_m\}$ is an orthonormal basis of $\Delta_n$, and $\hat{\rho}(E_a) \, \hat{s}_\alpha = \sum_{\beta=1}^m \omega_{a\alpha\beta} \, \hat{s}_\beta$ for $a=1,\hdots,n$. Consequently, $\hat{s} = (\hat{s}_1,\hdots,\hat{s}_m)$ is an $m$-tuple of harmonic spinors on $\Omega$ which satisfies the boundary condition $\chi \hat{s} = \hat{s}$. Thus, the kernel has dimension at least $1$.

We claim that the cokernel of the operator 
\[H^1(\Omega,\mathcal{E}) \to L^2(\Omega,\mathcal{E}) \oplus H^{\frac{1}{2}}(\Sigma,\mathcal{F}^-), \quad s \mapsto (\mathcal{D} s,s-\chi s)\] 
has dimension $0$. Suppose that this is false. In view of Proposition \ref{fredholm}, we can find a non-zero pair $(s,t)$ such that $s \in C^\infty(\Omega,\mathcal{E})$, $t \in C^\infty(\Sigma,\mathcal{F}^-)$, and  
\[\int_\Omega \langle s,\mathcal{D} u \rangle + \int_\Sigma \langle t,u-\chi u \rangle = 0\] 
for all $u \in C^\infty(\Omega,\mathcal{E})$. Integration by parts gives 
\[\int_\Omega \langle \mathcal{D} s,u \rangle - \int_\Sigma \langle \nu \cdot s,u \rangle + \int_\Sigma \langle t,u-\chi u \rangle = 0\] 
for all $u \in C^\infty(\Omega,\mathcal{E})$. From this, we deduce that $\mathcal{D} s = 0$ at each point in $\Omega$ and $\nu \cdot s = 2t$ at each point on $\Sigma$. Since $t$ takes values in $\mathcal{F}^-$, we conclude that the restriction $s|_\Sigma$ takes values in $\mathcal{F}^-$. In other words, $\chi s = -s$ at each point on $\Sigma$. Since $g$ is the Euclidean metric and $N$ is the Euclidean Gauss map of $\Sigma$, we have $H = \|dN\|_{\text{\rm tr}}$ at each point on $\Sigma$. Hence, Corollary \ref{rigidity.2} implies that $s = (s_1,\hdots,s_m)$ is parallel. Let us write $s_\alpha = \sum_{\beta=1}^m z_{\alpha\beta} \, \hat{s}_\beta$, where $z \in \text{\rm End}(\mathbb{C}^m)$ is a constant matrix. Since $\chi s = -s$ at each point on $\Sigma$, the matrix $z \in \text{\rm End}(\mathbb{C}^m)$ anti-commutes with the matrix $\sum_{a=1}^n \langle N(x),E_a \rangle \, \omega_a \in \text{\rm End}(\mathbb{C}^m)$ for each $x \in \Sigma$. Since the Gauss map $N: \Sigma \to S^{n-1}$ is surjective, the matrix $z \in \text{\rm End}(\mathbb{C}^m)$ anti-commutes with $\omega_a \in \text{\rm End}(\mathbb{C}^m)$ for each $a=1,\hdots,n$. Since $n$ is odd, Lemma \ref{algebraic.fact} implies that $z=0$. Thus, we conclude that $s=0$ at each point in $\Omega$. Since $\nu \cdot s = 2t$ at each point on $\Sigma$, it follows that $t=0$ at each point on $\Sigma$. This is a contradiction. Thus, the cokernel has dimension $0$. 

To summarize, if $g$ is the Euclidean metric and $N$ is the Euclidean Gauss map of $\Sigma$, then the index is at least $1$. 

We now turn to the general case. Let $g$ be an arbitrary Riemannian metric which is defined on an open set containing $\Omega$, and let $N: \Sigma \to S^{n-1}$ be a smooth map which is homotopic to the Euclidean Gauss map of $\Sigma$. Using the deformation invariance of the Fredholm index (cf. \cite{Hormander2}, Theorem 20.1.8), we conclude that the index is at least $1$. Note that in H\"ormander's setting, the vector bundles are fixed, whereas in our setting the vector bundle $\mathcal{F}^-$ depends on $g$ and $N$. To apply H\"ormander's results, we construct a bundle isomorphism from $\mathcal{F}^-$ to some fixed bundle. This completes the proof of Proposition \ref{index}. \\

\begin{remark} 
Under the assumptions of Proposition \ref{index}, we can show that the Fredholm index is equal to $1$. We will not need this stronger statement here.
\end{remark}

\section{Approximating a compact, convex polytope by smooth domains}

\label{smoothing.of.polytope}

Throughout this section, we assume that $n \geq 3$ is an integer and $\Omega$ is a compact, convex polytope in $\mathbb{R}^n$ with non-empty interior. We write $\Omega = \bigcap_{i \in I} \{u_i \leq 0\}$, where $u_i$, $i \in I$, is a finite collection of non-constant linear functions defined on $\mathbb{R}^n$. After eliminating redundant inequalities, we may assume that the following condition is satisfied. 

\begin{assumption}
\label{no.redundant.inequalities}
For each $i_0 \in I$, the set 
\[\{u_{i_0} > 0\} \cap \bigcap_{i \in I \setminus \{i_0\}} \{u_i \leq 0\}\] 
is non-empty. 
\end{assumption}

\begin{lemma}
\label{boundary.faces}
For each $i_0 \in I$, the set 
\[\{u_{i_0} = 0\} \cap \bigcap_{i \in I \setminus \{i_0\}} \{u_i < 0\}\] 
is non-empty. Moreover, this set is a dense subset of $\Omega \cap \{u_{i_0} = 0\}$.
\end{lemma}

\textbf{Proof.} 
In view of Assumption \ref{no.redundant.inequalities}, we can find a point $z_0 \in \mathbb{R}^n$ such that $u_{i_0}(z_0) > 0$ and $u_i(z_0) \leq 0$ for all $i \in I \setminus \{i_0\}$. Moreover, since $\Omega$ has non-empty interior, we can find a point $z_1 \in \mathbb{R}^n$ such that $u_i(z_1) < 0$ for all $i \in I$. We can find a real number $\tau \in (0,1)$ such that $(1-\tau) \, u_{i_0}(z_0) + \tau \, u_{i_0}(z_1) = 0$. Let $y := (1-\tau)z_0+\tau z_1$. Then 
\[y \in \{u_{i_0} = 0\} \cap \bigcap_{i \in I \setminus \{i_0\}} \{u_i < 0\}.\] 
This proves the first statement. To prove the second statement, we consider an arbitrary point $x \in \Omega \cap \{u_{i_0} = 0\}$. Then 
\[(1-t)x + ty \in \{u_{i_0} = 0\} \cap \bigcap_{i \in I \setminus \{i_0\}} \{u_i < 0\}\] 
for each $t \in (0,1]$. This completes the proof of Lemma \ref{boundary.faces}. \\

For each $i \in I$, we denote by $N_i \in S^{n-1}$ the outward-pointing unit normal vector to the halfspace $\{u_i \leq 0\}$ with respect to the Euclidean metric.

\begin{lemma}
\label{span}
We have $\mathbb{R}^n = \text{\rm span} \{N_i: i \in I\}$.
\end{lemma} 

\textbf{Proof.} 
Suppose that the assertion is false. Then we can find a non-zero vector in $\mathbb{R}^n$ which is orthogonal to $N_i$ for all $i \in I$. This implies that $\Omega$ is invariant under translations along that vector. This contradicts our assumption that $\Omega$ is compact. This completes the proof of Lemma \ref{span}. \\

Since $\Omega$ has non-empty interior, we can find a real number $\lambda_0 > 0$ such that 
\[\bigcap_{i \in I} \{u_i \leq -\lambda_0^{-1} \log |I|\} \neq \emptyset.\] 
For each $\lambda > \lambda_0$, we define 
\[\Omega_\lambda = \bigg \{ \sum_{i \in I} e^{\lambda u_i} \leq 1 \bigg \}.\] 
Clearly, 
\[\bigcap_{i \in I} \{u_i \leq -\lambda^{-1} \log |I|\} \subset \Omega_\lambda \subset \bigcap_{i \in I} \{u_i < 0\}\] 
for each $\lambda > \lambda_0$. In particular, 
\[\bigcup_{\lambda > \lambda_0} \Omega_\lambda = \bigcap_{i \in I} \{u_i < 0\} = \Omega \setminus \partial \Omega.\] 

\begin{lemma}
\label{Omega_lambda.is.smooth}
For each $\lambda > \lambda_0$, $\Omega_\lambda$ is a compact, convex domain in $\mathbb{R}^n$ with smooth boundary $\Sigma_\lambda = \partial \Omega_\lambda$.
\end{lemma}

\textbf{Proof.} 
Let us fix a real number $\lambda > \lambda_0$. It follows from Lemma \ref{span} that the function $\sum_{i \in I} e^{\lambda u_i}$ is strictly convex with respect to the Euclidean metric. Moreover, $\inf_{\partial \Omega} \sum_{i \in I} e^{\lambda u_i} > 1$. On the other hand, $\inf_\Omega \sum_{i \in I} e^{\lambda u_i} < 1$ since $\lambda > \lambda_0$. Consequently, we can find a point in the interior of $\Omega$ where the function $\sum_{i \in I} e^{\lambda u_i}$ attains its global minimum. From this, the assertion follows easily. This completes the proof of Lemma \ref{Omega_lambda.is.smooth}. \\

Let $g$ be a Riemannian metric which is defined on an open set containing $\Omega$. For each $i \in I$, $\nabla u_i$ will denote the gradient of $u_i$ with respect to the metric $g$; $D^2 u_i$ will denote the Hessian of $u_i$ with respect to the metric $g$; $|\nabla u_i|$ will denote the norm of the gradient of $u_i$ with respect to the metric $g$; and $\nu_i = \frac{\nabla u_i}{|\nabla u_i|}$ will denote the unit normal vector field, with respect to the metric $g$, to the level sets of $u_i$. 

\begin{lemma} 
\label{lower.bound.for.denominator.in.nu}
If $\lambda$ is sufficiently large, then $\inf_{\Sigma_\lambda} \big | \sum_{i \in I} e^{\lambda u_i} \, du_i \big | \geq C^{-1}$ for some large constant $C$ which is independent of $\lambda$.
\end{lemma}

\textbf{Proof.} 
We argue by contradiction. Suppose that the assertion is false. Then there exists a sequence of positive real numbers $\lambda_l \to \infty$ and a sequence of points $x_l \in \Sigma_{\lambda_l}$ such that $\big | \sum_{i \in I} e^{\lambda_l u_i} \, du_i \big | \leq l^{-1}$ at the point $x_l$. After passing to a subsequence, we may assume that the sequence $x_l$ converges to a point $x_0 \in \Omega$. Moreover, we may assume that, for each $i \in I$, the sequence $e^{\lambda_l u_i(x_l)}$ converges to a nonnegative real number $z_i$. Since $\sum_{i \in I} e^{\lambda_l u_i(x_l)} = 1$ for each $l$, we know that $\sum_{i \in I} z_i = 1$. Let $I_0 := \{i \in I: z_i > 0\}$. Clearly, $I_0$ is non-empty, and $u_i(x_0) = 0$ for all $i \in I_0$. Moreover, $\sum_{i \in I_0} z_i \, du_i = 0$ at the point $x_0$. On the other hand, since $\Omega$ is a convex set with non-empty interior, we can find a tangent vector $\xi \in T_{x_0} \Omega$ such that $du_i(\xi) > 0$ for all $i \in I_0$. This is a contradiction. This completes the proof of Lemma \ref{lower.bound.for.denominator.in.nu}. \\

\begin{lemma}
\label{lower.bound.for.denominator.in.N}
If $\lambda$ is sufficiently large, then $\inf_{\Sigma_\lambda} \big | \sum_{i \in I} e^{\lambda u_i} \, |\nabla u_i| \, N_i \big | \geq C^{-1}$ for some large constant $C$ which is independent of $\lambda$.
\end{lemma}

\textbf{Proof.} 
We argue by contradiction. Suppose that the assertion is false. Then there exists a sequence of positive real numbers $\lambda_l \to \infty$ and a sequence of points $x_l \in \Sigma_{\lambda_l}$ such that $\big | \sum_{i \in I} e^{\lambda_l u_i} \, |\nabla u_i| \, N_i \big | \leq l^{-1}$ at the point $x_l$. After passing to a subsequence, we may assume that the sequence $x_l$ converges to a point $x_0 \in \Omega$. Moreover, we may assume that, for each $i \in I$, the sequence $e^{\lambda_l u_i(x_l)} \, |\nabla u_i(x_l)|$ converges to a nonnegative real number $z_i$. Since $\sum_{i \in I} e^{\lambda_l u_i(x_l)} = 1$ for each $l$, we know that $\sum_{i \in I} z_i > 0$. Let $I_0 := \{i \in I: z_i > 0\}$. Clearly, $I_0$ is non-empty, and $u_i(x_0) = 0$ for all $i \in I_0$. Moreover, $\sum_{i \in I_0} z_i N_i = 0$ at the point $x_0$. On the other hand, since $\Omega$ is a convex set with non-empty interior, we can find a vector $\xi \in \mathbb{R}^n$ such that $\langle N_i,\xi \rangle > 0$ for all $i \in I_0$. This is a contradiction. This completes the proof of Lemma \ref{lower.bound.for.denominator.in.N}. \\

In the following, we assume that $\lambda$ is chosen sufficiently large so that the conclusions of Lemma \ref{lower.bound.for.denominator.in.nu} and Lemma \ref{lower.bound.for.denominator.in.N} hold. The outward-pointing unit normal vector to the domain $\Omega_\lambda$ with respect to the metric $g$ is given by 
\[\nu = \frac{\sum_{i \in I} e^{\lambda u_i} \, \nabla u_i}{\big | \sum_{i \in I} e^{\lambda u_i} \, \nabla u_i \big |} = \frac{\sum_{i \in I} e^{\lambda u_i} \, |\nabla u_i| \, \nu_i}{\big | \sum_{i \in I} e^{\lambda u_i} \, |\nabla u_i| \, \nu_i \big |}.\] 
This motivates the following definition: 

\begin{definition}
\label{definition.of.N}
We define a map $N: \Sigma_\lambda \to S^{n-1}$ by  
\[N = \frac{\sum_{i \in I} e^{\lambda u_i} \, |\nabla u_i| \, N_i}{\big | \sum_{i \in I} e^{\lambda u_i} \, |\nabla u_i| \, N_i \big |}.\] 
Recall that $|\nabla u_i|$ is computed with respect to the metric $g$. In particular, the map $N$ depends on the choice of the metric $g$.
\end{definition}

\begin{lemma}
\label{homotopy.class.of.N}
The map $N: \Sigma_\lambda \to S^{n-1}$ is homotopic to the Euclidean Gauss map of $\Sigma_\lambda$. 
\end{lemma} 

\textbf{Proof.} 
In the special case when $g$ is the Euclidean metric, the map $N: \Sigma_\lambda \to S^{n-1}$ coincides with the Gauss map of $\Sigma_\lambda$, and the assertion is trivial. To prove the assertion in general, we deform the metric $g$ to the Euclidean metric. \\

\begin{proposition} 
\label{formula.for.V}
Consider a point $x \in \Sigma_\lambda$. Let $\pi: T_x \Omega \to T_x \Omega$ denote the orthogonal projection to the orthogonal complement of $\nu$ with respect to $g$, and let $P: \mathbb{R}^n \to \mathbb{R}^n$ denote the orthogonal projection to the orthogonal complement of $N$ with respect to the Euclidean metric. Then $H - \|dN\|_{\text{\rm tr}} \geq V_\lambda$, where $H$ denotes the mean curvature of $\Sigma_\lambda$ with respect to the metric $g$ and the function $V_\lambda: \Sigma_\lambda \to \mathbb{R}$ is defined by 
\begin{align*} 
V_\lambda &= \lambda \, \frac{\sum_{i \in I} e^{\lambda u_i} \, |\nabla u_i|^2 \, |\pi(\nu_i)|^2}{\big | \sum_{i \in I} e^{\lambda u_i} \, |\nabla u_i| \, \nu_i \big |} - \lambda \, \frac{\sum_{i \in I} e^{\lambda u_i} \, |\nabla u_i|^2 \, |\pi(\nu_i)| \, |P(N_i)|}{\big | \sum_{i \in I} e^{\lambda u_i} \, |\nabla u_i| \, N_i \big |} \\ 
&+ \frac{\sum_{i \in I} e^{\lambda u_i} \, (\Delta u_i - (D^2 u_i)(\nu,\nu))}{\big | \sum_{i \in I} e^{\lambda u_i} \, |\nabla u_i| \, \nu_i \big |} - \frac{\sum_{i \in I} e^{\lambda u_i} \, |\nabla(|\nabla u_i|)| \, |P(N_i)|}{\big | \sum_{i \in I} e^{\lambda u_i} \, |\nabla u_i| \, N_i \big |}.
\end{align*} 
\end{proposition}

\textbf{Proof.} Let $\{e_1,\hdots,e_{n-1}\}$ denote a local orthonormal frame on $\Sigma_\lambda$ with respect to the metric $g$. The mean curvature of $\Sigma_\lambda$ with respect to $g$ is given by 
\begin{align*} 
H 
&= \lambda \, \frac{\sum_{k=1}^{n-1} \sum_{i \in I} e^{\lambda u_i} \, \langle \nabla u_i,e_k \rangle^2}{\big | \sum_{i \in I} e^{\lambda u_i} \, \nabla u_i \big |} + \frac{\sum_{k=1}^{n-1} \sum_{i \in I} e^{\lambda u_i} \, (D^2 u_i)(e_k,e_k)}{\big | \sum_{i \in I} e^{\lambda u_i} \, \nabla u_i \big |} \\ 
&= \lambda \, \frac{\sum_{i \in I} e^{\lambda u_i} \, |\pi(\nabla u_i)|^2}{\big | \sum_{i \in I} e^{\lambda u_i} \, \nabla u_i \big |} + \frac{\sum_{i \in I} e^{\lambda u_i} \, (\Delta u_i - (D^2 u_i)(\nu,\nu))}{\big | \sum_{i \in I} e^{\lambda u_i} \, \nabla u_i \big |} \\ 
&= \lambda \, \frac{\sum_{i \in I} e^{\lambda u_i} \, |\nabla u_i|^2 \, |\pi(\nu_i)|^2}{\big | \sum_{i \in I} e^{\lambda u_i} \, |\nabla u_i| \, \nu_i \big |} + \frac{\sum_{i \in I} e^{\lambda u_i} \, (\Delta u_i - (D^2 u_i)(\nu,\nu))}{\big | \sum_{i \in I} e^{\lambda u_i} \, |\nabla u_i| \, \nu_i \big |}. 
\end{align*} 
If $\xi$ is a tangent vector to $\Sigma_\lambda$, then  
\begin{align*} 
&dN(\xi) \\ 
&= \lambda \, \frac{\sum_{i \in I} e^{\lambda u_i} \, |\nabla u_i| \, \langle \nabla u_i,\xi \rangle \, P(N_i)}{\big | \sum_{i \in I} e^{\lambda u_i} \, |\nabla u_i| \, N_i \big |} + \frac{\sum_{i \in I} e^{\lambda u_i} \, \langle \nabla(|\nabla u_i|),\xi \rangle \, P(N_i)}{\big | \sum_{i \in I} e^{\lambda u_i} \, |\nabla u_i| \, N_i \big |} \\  
&= \lambda \, \frac{\sum_{i \in I} e^{\lambda u_i} \, |\nabla u_i|^2 \, \langle \pi(\nu_i),\xi \rangle \, P(N_i)}{\big | \sum_{i \in I} e^{\lambda u_i} \, |\nabla u_i| \, N_i \big |} + \frac{\sum_{i \in I} e^{\lambda u_i} \, \langle \nabla(|\nabla u_i|),\xi \rangle \, P(N_i)}{\big | \sum_{i \in I} e^{\lambda u_i} \, |\nabla u_i| \, N_i \big |}.
\end{align*} 
The trace norm of a linear map of the form $\xi \mapsto \langle X,\xi \rangle \, Y$ is given by $|X| \, |Y|$. Since the trace norm satisfies the triangle inequality, it follows that 
\begin{align*} 
&\|dN\|_{\text{\rm tr}} \\ 
&\leq \lambda \, \frac{\sum_{i \in I} e^{\lambda u_i} \, |\nabla u_i|^2 \, |\pi(\nu_i)| \, |P(N_i)|}{\big | \sum_{i \in I} e^{\lambda u_i} \, |\nabla u_i| \, N_i \big |} + \frac{\sum_{i \in I} e^{\lambda u_i} \, |\nabla(|\nabla u_i|)| \, |P(N_i)|}{\big | \sum_{i \in I} e^{\lambda u_i} \, |\nabla u_i| \, N_i \big |}.
\end{align*} 
Putting these facts together, the assertion follows. \\

\begin{proposition}
\label{negative.part.of.V.small.scale}
Suppose that the Matching Angle Hypothesis is satisfied. Then $\sup_{\Sigma_\lambda} \max \{-V_\lambda,0\} \leq o(\lambda)$ as $\lambda \to \infty$.
\end{proposition}

\textbf{Proof.} 
We argue by contradiction. Suppose that the assertion is false. Then there exists a sequence of positive real numbers $\lambda_l \to \infty$ and a sequence of points $x_l \in \Sigma_{\lambda_l}$ such that $\limsup_{l \to \infty} \lambda_l^{-1} \, V_{\lambda_l}(x_l) < 0$. After passing to a subsequence, we may assume that the sequence $x_l$ converges to a point $x_0 \in \Omega$. Moreover, we may assume that, for each $i \in I$, the sequence $e^{\lambda_l u_i(x_l)} \, |\nabla u_i(x_l)|$ converges to a nonnegative real number $z_i$. Since $\sum_{i \in I} e^{\lambda_l u_i(x_l)} = 1$ for each $l$, we know that $\sum_{i \in I} z_i > 0$. Let $I_0 := \{i \in I: z_i > 0\}$. Clearly, $I_0$ is non-empty, and $u_i(x_0) = 0$ for all $i \in I_0$. It follows from Lemma \ref{lower.bound.for.denominator.in.nu} and Lemma \ref{lower.bound.for.denominator.in.N} that $\big | \sum_{i \in I_0} z_i \nu_i \big | > 0$ and $\big | \sum_{i \in I_0} z_i N_i \big | > 0$ at the point $x_0$.

We now invoke the Matching Angle Hypothesis. Hence, for all $i_1,i_2 \in I_0$, we have $\langle \nu_{i_1},\nu_{i_2} \rangle = \langle N_{i_1},N_{i_2} \rangle$ at the point $x_0$. Let $\pi: T_{x_0} \Omega \to T_{x_0} \Omega$ denote the orthogonal projection to the orthogonal complement of $\sum_{i \in I_0} z_i \nu_i$, and let $P: \mathbb{R}^n \to \mathbb{R}^n$ denote the orthogonal projection to the orthogonal complement of $\sum_{i \in I_0} z_i N_i$. For each $j \in I_0$, we have 
\begin{align*} 
\Big | \sum_{i \in I_0} z_i N_i \Big |^2 \, |P(N_j)|^2 
&= \Big | \sum_{i \in I_0} z_i N_i \Big |^2 - \Big \langle \sum_{i \in I_0} z_i N_i,N_j \Big \rangle^2 \\ 
&= \Big | \sum_{i \in I_0} z_i \nu_i \Big |^2 - \Big \langle \sum_{i \in I_0} z_i \nu_i,\nu_j \Big \rangle^2 \\ 
&= \Big | \sum_{i \in I_0} z_i \nu_i \Big |^2 \, |\pi(\nu_j)|^2 
\end{align*}
at the point $x_0$. Moreover, 
\[\Big | \sum_{i \in I_0} z_i N_i \Big |^2 = \Big | \sum_{i \in I_0} z_i \nu_i \Big |^2.\] 
at the point $x_0$. Consequently, for each $j \in I_0$, we obtain  
\[\frac{|\pi(\nu_j)| \, |P(N_j)|}{\big | \sum_{i \in I_0} z_i N_i \big |} = \frac{|\pi(\nu_j)|^2}{\big | \sum_{i \in I_0} z_i \nu_i \big |}\] 
at the point $x_0$. This implies 
\[\frac{\sum_{i \in I_0} z_i \, |\nabla u_i| \, |\pi(\nu_i)| \, |P(N_i)|}{\big | \sum_{i \in I_0} z_i N_i \big |} = \frac{\sum_{i \in I_0} z_i \, |\nabla u_i| \, |\pi(\nu_i)|^2}{\big | \sum_{i \in I_0} z_i \nu_i \big |}\] 
at the point $x_0$. Using Proposition \ref{formula.for.V}, we conclude that $\lambda_l^{-1} \, V_{\lambda_l}(x_l) \to 0$ as $l \to \infty$. This is a contradiction. \\

In the remainder of this section, we estimate the $L^\sigma$-norm of $\max \{-V_\lambda,0\}$ on $\Sigma_\lambda \cap B_r(p)$, where $\sigma \in [1,\frac{3}{2})$ is a fixed exponent and $B_r(p)$ denotes a Euclidean ball of radius $r$ centered at a point $p \in \mathbb{R}^n$. We first recall a basic fact about the area of convex hypersurfaces in $\mathbb{R}^n$.

\begin{lemma} 
\label{area}
Let $B_r(p)$ denote a Euclidean ball of radius $r$. Then the intersection $\Sigma_\lambda \cap B_r(p)$ has area at most $C r^{n-1}$, where $C$ depends only on $n$.
\end{lemma}

\textbf{Proof.} 
The hypersurface $\Sigma_\lambda$ bounds a convex domain in Euclidean space. This implies that $\Sigma_\lambda$ is outward-minimizing with respect to the Euclidean metric. From this, the assertion follows. \\

\begin{definition}
Consider three pairwise distinct elements $i_1,i_2,i_3 \in I$. We denote by $G_\lambda^{(i_1,i_2,i_3)}$ the set of all points $x \in \Sigma_\lambda$ with the property that $u_{i_1}(x) \geq u_{i_2}(x)\geq u_{i_3}(x)$ and $u_{i_3}(x) \geq u_i(x)$ for each $i \in I \setminus \{i_1,i_2,i_3\}$. 
\end{definition}

Clearly, $\Sigma_\lambda = \bigcup_{i_1,i_2,i_3} G_\lambda^{(i_1,i_2,i_3)}$, where the union is taken over all triplets $(i_1,i_2,i_3) \in I \times I \times I$ such that $i_1,i_2,i_3$ are pairwise distinct. Given three pairwise distinct elements $i_1,i_2,i_3 \in I$, we shall estimate the $L^\sigma$-norm of $\max \{-V_\lambda,0\}$ on the set $G_\lambda^{(i_1,i_2,i_3)} \cap B_r(p)$. To that end, we decompose the set $G_\lambda^{(i_1,i_2,i_3)}$ into three subsets. Roughly speaking, the first subset consists of points that are close to one of the $(n-1)$-dimensional boundary faces of $\Omega$, but stay away from the $(n-2)$-dimensional edges of $\Omega$. The second subset consists of points that are close to one of the $(n-2)$-dimensional edges of $\Omega$, but stay away from the $(n-3)$-dimensional corners of $\Omega$. The third set consists of points that are close to one of the $(n-3)$-dimensional corners of $\Omega$. 

\begin{lemma}
\label{step.1}
For each $i \in I$, we assume that the mean curvature of the hypersurface $\{u_i=0\}$ with respect to $g$ is nonnegative at each point in $\Omega \cap \{u_i=0\}$. Let us fix an exponent $\sigma \in [1,\frac{3}{2})$, and let $B_r(p)$ denote a Euclidean ball of radius $0 < r \leq 1$. If $\lambda r$ is sufficiently large, then 
\[\bigg ( r^{\sigma+1-n} \int_{G_\lambda^{(i_1,i_2,i_3)} \cap \{u_{i_2} \leq -\lambda^{-\frac{7}{8}} r^{\frac{1}{8}}\} \cap B_r(p)} (\max \{-V_\lambda,0\})^\sigma \bigg )^{\frac{1}{\sigma}} \leq C\lambda r \, e^{-(\lambda r)^{\frac{1}{8}}}\] 
for all pairwise distinct elements $i_1,i_2,i_3 \in I$. The constant $C$ may depend on $\Omega$ and $g$, but not on $\lambda$. 
\end{lemma}

\textbf{Proof.} 
Let us consider an arbitrary point $x \in G_\lambda^{(i_1,i_2,i_3)}$ with $u_{i_2}(x) \leq -\lambda^{-\frac{7}{8}} r^{\frac{1}{8}}$. By definition of $G_\lambda^{(i_1,i_2,i_3)}$, it follows that $u_i(x) \leq -\lambda^{-\frac{7}{8}} r^{\frac{1}{8}}$ for all $i \in I \setminus \{i_1\}$. Using the identity $\sum_{i \in I} e^{\lambda u_i(x)} = 1$, we obtain $e^{\lambda u_{i_1}(x)} \geq 1 - C \, e^{-(\lambda r)^{\frac{1}{8}}}$, hence $u_{i_1}(x) \geq -C\lambda^{-1} \, e^{-(\lambda r)^{\frac{1}{8}}}$. Moreover, $|\nu-\nu_{i_1}| \leq C \, e^{-(\lambda r)^{\frac{1}{8}}}$ and $|N-N_{i_1}| \leq C \, e^{-(\lambda r)^{\frac{1}{8}}}$ at the point $x$. From this, we deduce that $|\pi(\nu_{i_1})| = |\pi(\nu_{i_1} - \nu)| \leq C \, e^{-(\lambda r)^{\frac{1}{8}}}$ and $|P(N_{i_1})| = |P(N_{i_1} - N)| \leq C \, e^{-(\lambda r)^{\frac{1}{8}}}$ at the point $x$. This gives 
\begin{equation} 
\label{term.1}
\sum_{i \in I} e^{\lambda u_i} \, |\nabla u_i|^2 \, |\pi(\nu_i)| \, |P(N_i)| \leq C \, e^{-(\lambda r)^{\frac{1}{8}}} 
\end{equation}
and 
\begin{equation} 
\label{term.2}
\sum_{i \in I} e^{\lambda u_i} \, |\nabla(|\nabla u_i|)| \, |P(N_i)| \leq C \, e^{-(\lambda r)^{\frac{1}{8}}} 
\end{equation}
at the point $x$. Moreover, 
\begin{align} 
\label{term.3}
&\sum_{i \in I} e^{\lambda u_i} \, (\Delta u_i - (D^2 u_i)(\nu,\nu)) \notag \\
&\geq e^{\lambda u_{i_1}} \, (\Delta u_{i_1} - (D^2 u_{i_1})(\nu_{i_1},\nu_{i_1})) - C \, e^{-(\lambda r)^{\frac{1}{8}}} 
\end{align}
at the point $x$. Since $u_{i_1}(x) \geq -C\lambda^{-1} \, e^{-(\lambda r)^{\frac{1}{8}}}$, we can find a point $y \in \mathbb{R}^n$ such that $u_{i_1}(y) = 0$ and $d_{\text{\rm eucl}}(x,y) \leq C\lambda^{-1} \, e^{-(\lambda r)^{\frac{1}{8}}}$. This implies $u_i(y) \leq u_i(x) + C\lambda^{-1} \, e^{-(\lambda r)^{\frac{1}{8}}} \leq -\lambda^{-\frac{7}{8}} r^{\frac{1}{8}} + C\lambda^{-1} \, e^{-(\lambda r)^{\frac{1}{8}}}$ for all $i \in I \setminus \{i_1\}$. In particular, if $\lambda r$ is sufficiently large, then $u_i(y) \leq 0$ for all $i \in I \setminus \{i_1\}$. Thus, $y \in \Omega \cap \{u_{i_1}=0\}$. By assumption, the mean curvature of the hypersurface $\{u_{i_1}=0\}$ at the point $y$ is nonnegative. This implies 
\[\Delta u_{i_1} - (D^2 u_{i_1})(\nu_{i_1},\nu_{i_1}) \geq 0\] 
at the point $y$. Consequently, 
\begin{equation} 
\label{mean.curvature.of.level.set.of.u_i_1}
\Delta u_{i_1} - (D^2 u_{i_1})(\nu_{i_1},\nu_{i_1}) \geq -C\lambda^{-1} \, e^{-(\lambda r)^{\frac{1}{8}}} 
\end{equation}
at the point $x$. Combining (\ref{term.3}) and (\ref{mean.curvature.of.level.set.of.u_i_1}), we obtain 
\begin{equation} 
\label{term.3.new}
\sum_{i \in I} e^{\lambda u_i} \, (\Delta u_i - (D^2 u_i)(\nu,\nu)) \geq -C \, e^{-(\lambda r)^{\frac{1}{8}}}. 
\end{equation}
Using (\ref{term.1}), (\ref{term.2}), and (\ref{term.3.new}), we conclude that 
\[V_\lambda(x) \geq -C\lambda \, e^{-(\lambda r)^{\frac{1}{8}}}\] 
for each point $x \in G_\lambda^{(i_1,i_2,i_3)} \cap \{u_{i_2} \leq -\lambda^{-\frac{7}{8}} r^{\frac{1}{8}}\}$. By Lemma \ref{area}, the intersection $\Sigma_\lambda \cap B_r(p)$ has area at most $C r^{n-1}$. Consequently, 
\[\bigg ( r^{\sigma+1-n} \int_{G_\lambda^{(i_1,i_2,i_3)} \cap \{u_{i_2} \leq -\lambda^{-\frac{7}{8}} r^{\frac{1}{8}}\} \cap B_r(p)} (\max \{-V_\lambda,0\})^\sigma \bigg )^{\frac{1}{\sigma}} \leq C\lambda r \, e^{-(\lambda r)^{\frac{1}{8}}}.\] 
This completes the proof of Lemma \ref{step.1}. \\

\begin{lemma} 
\label{step.2}
Assume that the Matching Angle Hypothesis holds.  Let us fix an exponent $\sigma \in [1,\frac{3}{2})$, and let $B_r(p)$ denote a Euclidean ball of radius $0 < r \leq 1$. If $\lambda r$ is sufficiently large, then 
\begin{align*} 
&\bigg ( r^{\sigma+1-n} \int_{G_\lambda^{(i_1,i_2,i_3)} \cap \{u_{i_2} \geq -\lambda^{-\frac{7}{8}} r^{\frac{1}{8}}\} \cap \{u_{i_3} \leq -\lambda^{-\frac{3}{4}} r^{\frac{1}{4}}\} \cap B_r(p)} (\max \{-V_\lambda,0\})^\sigma \bigg )^{\frac{1}{\sigma}} \\ 
&\leq C \, (\lambda r)^{\frac{1}{8}-\frac{7}{8\sigma}} 
\end{align*}
for all pairwise distinct elements $i_1,i_2,i_3 \in I$. The constant $C$ may depend on $\Omega$ and $g$, but not on $\lambda$. 
\end{lemma}

\textbf{Proof.} 
We distinguish two cases: 

\textit{Case 1:} Suppose that $\Omega \cap \{u_{i_1}=0\} \cap \{u_{i_2}=0\} = \emptyset$. We can find a positive real number $\delta$ such that $\Omega \cap \{u_{i_1} \geq -\delta\} \cap \{u_{i_2} \geq -\delta\} = \emptyset$. If $\lambda r$ is sufficiently large, then $\lambda^{-\frac{7}{8}} r^{\frac{1}{8}} \leq (\lambda r)^{-\frac{7}{8}} \leq \delta$. This implies 
\begin{align*} 
&G_\lambda^{(i_1,i_2,i_3)} \cap \{u_{i_2} \geq -\lambda^{-\frac{7}{8}} r^{\frac{1}{8}}\} \\ 
&\subset \Sigma_\lambda \cap \{u_{i_1} \geq -\delta\} \cap \{u_{i_2} \geq -\delta\} = \emptyset. 
\end{align*} 
Hence, the assertion is trivially true in this case. 

\textit{Case 2:} Suppose that $\Omega \cap \{u_{i_1}=0\} \cap \{u_{i_2}=0\} \neq \emptyset$. It follows from Assumption \ref{no.redundant.inequalities} that the hypersurfaces $\{u_{i_1}=0\}$ and $\{u_{i_2}=0\}$ intersect transversally.

Let us consider an arbitrary point $x \in G_\lambda^{(i_1,i_2,i_3)}$ with $u_{i_2}(x) \geq -\lambda^{-\frac{7}{8}} r^{\frac{1}{8}}$ and $u_{i_3}(x) \leq -\lambda^{-\frac{3}{4}} r^{\frac{1}{4}}$. Clearly, $u_{i_1}(x) \geq -\lambda^{-\frac{7}{8}} r^{\frac{1}{8}}$ by definition of $G_\lambda^{(i_1,i_2,i_3)}$. By transversality, we can find a point $y \in \mathbb{R}^n$ such that $u_{i_1}(y) = u_{i_2}(y) = 0$ and $d_{\text{\rm eucl}}(x,y) \leq C \lambda^{-\frac{7}{8}} r^{\frac{1}{8}}$. This implies $u_i(y) \leq u_i(x) + C \lambda^{-\frac{7}{8}} r^{\frac{1}{8}} \leq -\lambda^{-\frac{3}{4}} r^{\frac{1}{4}} + C \lambda^{-\frac{7}{8}} r^{\frac{1}{8}}$ for all $i \in I \setminus \{i_1,i_2\}$. In particular, if $\lambda r$ is sufficiently large, then $u_i(y) \leq 0$ for all $i \in I \setminus \{i_1,i_2\}$. Thus, $y \in \Omega \cap \{u_{i_1}=0\} \cap \{u_{i_2}=0\}$. The Matching Angle Hypothesis implies that $\langle \nu_{i_1},\nu_{i_2} \rangle = \langle N_{i_1},N_{i_2} \rangle$ at the point $y$. Consequently, $|\langle \nu_{i_1},\nu_{i_2} \rangle - \langle N_{i_1},N_{i_2} \rangle| \leq C \lambda^{-\frac{7}{8}} r^{\frac{1}{8}} \leq C \, (\lambda r)^{-\frac{7}{8}}$ at the point $x$. For each $j \in \{i_1,i_2\}$, we have 
\begin{align*} 
&\Big | \sum_{i \in I} e^{\lambda u_i} \, |\nabla u_i| \, N_i \Big |^2 \, |P(N_j)|^2 \\ 
&= \Big | \sum_{i \in I} e^{\lambda u_i} \, |\nabla u_i| \, N_i \Big |^2 - \Big \langle \sum_{i \in I} e^{\lambda u_i} \, |\nabla u_i| \, N_i,N_j \Big \rangle^2 \\ 
&\leq \Big | \sum_{i \in I} e^{\lambda u_i} \, |\nabla u_i| \, \nu_i \Big |^2 - \Big \langle \sum_{i \in I} e^{\lambda u_i} \, |\nabla u_i| \, \nu_i,\nu_j \Big \rangle^2 + C \, (\lambda r)^{-\frac{7}{8}} \\ 
&= \Big | \sum_{i \in I} e^{\lambda u_i} \, |\nabla u_i| \, \nu_i \Big |^2 \, |\pi(\nu_j)|^2 + C \, (\lambda r)^{-\frac{7}{8}} 
\end{align*} 
at the point $x$. Hence, for each $j \in \{i_1,i_2\}$, we obtain 
\begin{align} 
\label{numerator}
&\Big | \sum_{i \in I} e^{\lambda u_i} \, |\nabla u_i| \, \nu_i \Big | \, \Big | \sum_{i \in I} e^{\lambda u_i} \, |\nabla u_i| \, N_i \Big | \, |\pi(\nu_j)| \, |P(N_j)| \notag \\ 
&\leq \frac{1}{2} \, \Big | \sum_{i \in I} e^{\lambda u_i} \, |\nabla u_i| \, \nu_i \Big |^2 \, |\pi(\nu_j)|^2 + \frac{1}{2} \, \Big | \sum_{i \in I} e^{\lambda u_i} \, |\nabla u_i| \, N_i \Big |^2 \, |P(N_j)|^2 \\ 
&\leq \Big | \sum_{i \in I} e^{\lambda u_i} \, |\nabla u_i| \, \nu_i \Big |^2 \, |\pi(\nu_j)|^2 +  C \, (\lambda r)^{-\frac{7}{8}} \notag 
\end{align}
at the point $x$. Moreover, 
\begin{equation} 
\label{denominator}
\Big | \sum_{i \in I} e^{\lambda u_i} \, |\nabla u_i| \, N_i \Big |^2 \geq \Big | \sum_{i \in I} e^{\lambda u_i} \, |\nabla u_i| \, \nu_i \Big |^2 - C \, (\lambda r)^{-\frac{7}{8}} 
\end{equation}
at the point $x$. In the next step, we divide the inequality (\ref{numerator}) by (\ref{denominator}). It follows from Lemma \ref{lower.bound.for.denominator.in.nu} that $\big | \sum_{i \in I} e^{\lambda u_i} \, |\nabla u_i| \, \nu_i \big | \geq C^{-1}$ at the point $x$. Lemma \ref{lower.bound.for.denominator.in.N} implies that $\big | \sum_{i \in I} e^{\lambda u_i} \, |\nabla u_i| \, N_i \big | \geq C^{-1}$ at the point $x$. Consequently, for each $j \in \{i_1,i_2\}$, we have 
\[\frac{|\pi(\nu_j)| \, |P(N_j)|}{\big | \sum_{i \in I} e^{\lambda u_i} \, |\nabla u_i| \, N_i \big |} \leq \frac{|\pi(\nu_j)|^2}{\big | \sum_{i \in I} e^{\lambda u_i} \, |\nabla u_i| \, \nu_i \big |} + C \, (\lambda r)^{-\frac{7}{8}}\]
at the point $x$. This implies 
\[\frac{\sum_{i \in I} e^{\lambda u_i} \, |\nabla u_i|^2 \, |\pi(\nu_i)| \, |P(N_i)|}{\big | \sum_{i \in I} e^{\lambda u_i} \, |\nabla u_i| \, N_i \big |} \leq \frac{\sum_{i \in I} e^{\lambda u_i} \, |\nabla u_i|^2 \, |\pi(\nu_i)|^2}{\big | \sum_{i \in I} e^{\lambda u_i} \, |\nabla u_i| \, \nu_i \big |} + C \, (\lambda r)^{-\frac{7}{8}}\]
at the point $x$. Thus, we conclude that 
\[V_\lambda(x) \geq -C \lambda^{\frac{1}{8}} r^{-\frac{7}{8}}\] 
for each point $x \in G_\lambda^{(i_1,i_2,i_3)} \cap \{u_{i_2} \geq -\lambda^{-\frac{7}{8}} r^{\frac{1}{8}}\} \cap \{u_{i_3} \leq -\lambda^{-\frac{3}{4}} r^{\frac{1}{4}}\}$. By transversality, the set $\{0 \geq u_{i_1} \geq -\lambda^{-\frac{7}{8}} r^{\frac{1}{8}}\} \cap \{0 \geq u_{i_2} \geq -\lambda^{-\frac{7}{8}} r^{\frac{1}{8}}\} \cap B_r(p)$ can be covered by $C \, (\lambda r)^{\frac{7(n-2)}{8}}$ Euclidean balls of radius $\lambda^{-\frac{7}{8}} r^{\frac{1}{8}}$. By Lemma \ref{area}, the intersection of $\Sigma_\lambda$ with each ball of radius $\lambda^{-\frac{7}{8}} r^{\frac{1}{8}}$ has area at most $C \, (\lambda r)^{-\frac{7(n-1)}{8}} \, r^{n-1}$. This implies that the set $\Sigma_\lambda \cap \{u_{i_1} \geq -\lambda^{-\frac{7}{8}} r^{\frac{1}{8}}\} \cap \{u_{i_2} \geq -\lambda^{-\frac{7}{8}} r^{\frac{1}{8}}\} \cap B_r(p)$ has area at most $C \, (\lambda r)^{-\frac{7}{8}} \, r^{n-1}$. Since 
\begin{align*} 
&G_\lambda^{(i_1,i_2,i_3)} \cap \{u_{i_2} \geq -\lambda^{-\frac{7}{8}} r^{\frac{1}{8}}\} \cap B_r(p) \\ 
&\subset \Sigma_\lambda \cap \{u_{i_1} \geq -\lambda^{-\frac{7}{8}} r^{\frac{1}{8}}\} \cap \{u_{i_2} \geq -\lambda^{-\frac{7}{8}} r^{\frac{1}{8}}\} \cap B_r(p), 
\end{align*} 
it follows that 
\begin{align*} 
&\bigg ( r^{\sigma+1-n} \int_{G_\lambda^{(i_1,i_2,i_3)} \cap \{u_{i_2} \geq -\lambda^{-\frac{7}{8}} r^{\frac{1}{8}}\} \cap \{u_{i_3} \leq -\lambda^{-\frac{3}{4}} r^{\frac{1}{4}}\} \cap B_r(p)} (\max \{-V_\lambda,0\})^\sigma \bigg )^{\frac{1}{\sigma}} \\ 
&\leq C \, (\lambda r)^{\frac{1}{8}-\frac{7}{8\sigma}}. 
\end{align*}
This completes the proof of Lemma \ref{step.2}. \\

\begin{lemma} 
\label{step.3}
Let us fix an exponent $\sigma \in [1,\frac{3}{2})$, and let $B_r(p)$ denote a Euclidean ball of radius $0 < r \leq 1$. If $\lambda r$ is sufficiently large, then 
\[\bigg ( r^{\sigma+1-n} \int_{G_\lambda^{(i_1,i_2,i_3)} \cap \{u_{i_3} \geq -\lambda^{-\frac{3}{4}} r^{\frac{1}{4}}\} \cap B_r(p)} (\max \{-V_\lambda,0\})^\sigma \bigg )^{\frac{1}{\sigma}} \leq C \, (\lambda r)^{1-\frac{3}{2\sigma}}\] 
for all pairwise distinct elements $i_1,i_2,i_3 \in I$. The constant $C$ may depend on $\Omega$ and $g$, but not on $\lambda$. 
\end{lemma}

\textbf{Proof.} 
We distinguish two cases: 

\textit{Case 1:} Suppose that $\Omega \cap \{u_{i_1}=0\} \cap \{u_{i_2}=0\} \cap \{u_{i_3}=0\} = \emptyset$. We can find a positive real number $\delta$ such that $\Omega \cap \{u_{i_1} \geq -\delta\} \cap \{u_{i_2} \geq -\delta\} \cap \{u_{i_3} \geq -\delta\} = \emptyset$. If $\lambda r$ is sufficiently large, then $\lambda^{-\frac{3}{4}} r^{\frac{1}{4}} \leq (\lambda r)^{-\frac{3}{4}} \leq \delta$. This implies 
\begin{align*} 
&G_\lambda^{(i_1,i_2,i_3)} \cap \{u_{i_3} \geq -\lambda^{-\frac{3}{4}} r^{\frac{1}{4}}\} \\ 
&\subset \Sigma_\lambda \cap \{u_{i_1} \geq -\delta\} \cap \{u_{i_2} \geq -\delta\} \cap \{u_{i_3} \geq -\delta\} = \emptyset. 
\end{align*} 
Hence, the assertion is trivially true in this case. 

\textit{Case 2:} Suppose that $\Omega \cap \{u_{i_1}=0\} \cap \{u_{i_2}=0\} \cap \{u_{i_3}=0\} \neq \emptyset$. It follows from Assumption \ref{no.redundant.inequalities} that the hypersurfaces $\{u_{i_1}=0\}$, $\{u_{i_2}=0\}$, and $\{u_{i_3}=0\}$ intersect transversally.

Let us consider an arbitrary point $x \in G_\lambda^{(i_1,i_2,i_3)}$ with $u_{i_3}(x) \geq -\lambda^{-\frac{3}{4}} r^{\frac{1}{4}}$. Clearly,  
\[V_\lambda(x) \geq -C\lambda\] 
for all points $x \in G_\lambda^{(i_1,i_2,i_3)} \cap \{u_{i_3} \geq -\lambda^{-\frac{3}{4}} r^{\frac{1}{4}}\}$. By transversality, the set $\{0 \geq u_{i_1} \geq -\lambda^{-\frac{3}{4}} r^{\frac{1}{4}}\} \cap \{0 \geq u_{i_2} \geq -\lambda^{-\frac{3}{4}} r^{\frac{1}{4}}\} \cap \{0 \geq u_{i_3} \geq -\lambda^{-\frac{3}{4}} r^{\frac{1}{4}}\} \cap B_r(p)$ can be covered by $C \, (\lambda r)^{\frac{3(n-3)}{4}}$ Euclidean balls of radius $\lambda^{-\frac{3}{4}} r^{\frac{1}{4}}$. By Lemma \ref{area}, the intersection of $\Sigma_\lambda$ with each ball of radius $\lambda^{-\frac{3}{4}} r^{\frac{1}{4}}$ has area at most $C \, (\lambda r)^{-\frac{3(n-1)}{4}} \, r^{n-1}$. This implies that the set $\Sigma_\lambda \cap \{u_{i_1} \geq -\lambda^{-\frac{3}{4}} r^{\frac{1}{4}}\} \cap \{u_{i_2} \geq -\lambda^{-\frac{3}{4}} r^{\frac{1}{4}}\} \cap \{u_{i_3} \geq -\lambda^{-\frac{3}{4}} r^{\frac{1}{4}}\} \cap B_r(p)$ has area at most $C \, (\lambda r)^{-\frac{3}{2}} \, r^{n-1}$. Since 
\begin{align*} 
&G_\lambda^{(i_1,i_2,i_3)} \cap \{u_{i_3} \geq -\lambda^{-\frac{3}{4}} r^{\frac{1}{4}}\} \cap B_r(p) \\ 
&\subset \Sigma_\lambda \cap \{u_{i_1} \geq -\lambda^{-\frac{3}{4}} r^{\frac{1}{4}}\} \cap \{u_{i_2} \geq -\lambda^{-\frac{3}{4}} r^{\frac{1}{4}}\} \cap \{u_{i_3} \geq -\lambda^{-\frac{3}{4}} r^{\frac{1}{4}}\} \cap B_r(p), 
\end{align*}
it follows that 
\[\bigg ( r^{\sigma+1-n} \int_{G_\lambda^{(i_1,i_2,i_3)} \cap \{u_{i_3} \geq -\lambda^{-\frac{3}{4}} r^{\frac{1}{4}}\} \cap B_r(p)} (\max \{-V_\lambda,0\})^\sigma \bigg )^{\frac{1}{\sigma}} \leq C \, (\lambda r)^{1-\frac{3}{2\sigma}}.\] 
This completes the proof of Lemma \ref{step.3}. \\

\begin{proposition}
\label{negative.part.of.V.large.scale}
For each $i \in I$, we assume that the mean curvature of the hypersurface $\{u_i=0\}$ with respect to $g$ is nonnegative at each point in $\Omega \cap \{u_i=0\}$. Moreover, we assume that the Matching Angle Hypothesis is satisfied. Let us fix an exponent $\sigma \in [1,\frac{3}{2})$, and let $B_r(p)$ denote a Euclidean ball of radius $0 < r \leq 1$. If $\lambda r$ is sufficiently large, then 
\begin{align*} 
&\bigg ( r^{\sigma+1-n} \int_{\Sigma_\lambda \cap B_r(p)} (\max \{-V_\lambda,0\})^\sigma \bigg )^{\frac{1}{\sigma}} \\ 
&\leq C\lambda r \, e^{-(\lambda r)^{\frac{1}{8}}} + C \, (\lambda r)^{\frac{1}{8}-\frac{7}{8\sigma}} + C \, (\lambda r)^{1-\frac{3}{2\sigma}}. 
\end{align*}
The constant $C$ may depend on $\Omega$ and $g$, but not on $\lambda$. 
\end{proposition}

\textbf{Proof.} 
Combining Lemma \ref{step.1}, Lemma \ref{step.2}, and Lemma \ref{step.3}, we conclude that 
\begin{align*} 
&\bigg ( r^{\sigma+1-n} \int_{G_\lambda^{(i_1,i_2,i_3)} \cap B_r(p)} (\max \{-V_\lambda,0\})^\sigma \bigg )^{\frac{1}{\sigma}} \\ 
&\leq C\lambda r \, e^{-(\lambda r)^{\frac{1}{8}}} + C \, (\lambda r)^{\frac{1}{8}-\frac{7}{8\sigma}} + C \, (\lambda r)^{1-\frac{3}{2\sigma}} 
\end{align*}
for all pairwise distinct elements $i_1,i_2,i_3 \in I$. On the other hand, $\Sigma_\lambda = \bigcup_{i_1,i_2,i_3} G_\lambda^{(i_1,i_2,i_3)}$, where the union is taken over all triplets $(i_1,i_2,i_3) \in I \times I \times I$ such that $i_1,i_2,i_3$ are pairwise distinct. Hence, the assertion follows by summation over $i_1,i_2,i_3$. This completes the proof of Proposition \ref{negative.part.of.V.large.scale}. \\

\begin{corollary}
\label{negative.part.of.V}
For each $i \in I$, we assume that the mean curvature of the hypersurface $\{u_i=0\}$ with respect to $g$ is nonnegative at each point in $\Omega \cap \{u_i=0\}$. Moreover, we assume that the Matching Angle Hypothesis is satisfied. Let us fix an exponent $\sigma \in [1,\frac{3}{2})$. Then 
\[\sup_{p \in \mathbb{R}^n} \sup_{0 < r \leq 1} \bigg ( r^{\sigma+1-n} \int_{\Sigma_\lambda \cap B_r(p)} (\max \{-V_\lambda,0\})^\sigma \bigg )^{\frac{1}{\sigma}} \to 0\] 
as $\lambda \to \infty$.
\end{corollary}

\textbf{Proof.} 
Let us consider an arbitrary sequence $\lambda_l \to \infty$. By Proposition \ref{negative.part.of.V.small.scale}, we can find a sequence of positive real numbers $\delta_l \to 0$ such that 
\[(\delta_l \lambda_l)^{-1} \sup_{\Sigma_{\lambda_l}} \max \{-V_{\lambda_l},0\} \to 0\] 
as $l \to \infty$. Using Lemma \ref{area}, we obtain 
\begin{align*} 
&\sup_{p \in \mathbb{R}^n} \sup_{0 < r \leq (\delta_l \lambda_l)^{-1}} \bigg ( r^{\sigma+1-n} \int_{\Sigma_{\lambda_l} \cap B_r(p)} (\max \{-V_{\lambda_l},0\})^\sigma \bigg )^{\frac{1}{\sigma}} \\ 
&\leq C \, (\delta_l \lambda_l)^{-1} \sup_{\Sigma_{\lambda_l}} \max \{-V_{\lambda_l},0\} \to 0 
\end{align*}
as $l \to \infty$. On the other hand, it follows from Proposition \ref{negative.part.of.V.large.scale} that 
\[\sup_{p \in \mathbb{R}^n} \sup_{(\delta_l \lambda_l)^{-1} \leq r \leq 1} \bigg ( r^{\sigma+1-n} \int_{\Sigma_{\lambda_l} \cap B_r(p)} (\max \{-V_{\lambda_l},0\})^\sigma \bigg )^{\frac{1}{\sigma}} \to 0\] 
as $l \to \infty$. Putting these facts together, the assertion follows. \\

\section{Proof of the Theorem \ref{main.theorem}}

It suffices to prove Theorem \ref{main.theorem} in the odd-dimensional case. (The even-dimensional case can be reduced to the odd-dimensional case by considering the Cartesian product $\Omega \times [-1,1] \subset \mathbb{R}^{n+1}$.) Suppose that $n \geq 3$ is an odd integer, and $\Omega$ is a compact, convex polytope in $\mathbb{R}^n$ with non-empty interior. We write $\Omega = \bigcap_{i \in I} \{u_i \leq 0\}$, where $u_i$, $i \in I$, is a finite collection of non-constant linear functions defined on $\mathbb{R}^n$. Let $g$ be a Riemannian metric which is defined on an open set containing $\Omega$ and has nonnegative scalar curvature at each point in $\Omega$. For each $i \in I$, we assume that the mean curvature of the hypersurface $\{u_i=0\}$ with respect to $g$ is nonnegative at each point in $\Omega \cap \{u_i=0\}$. Moreover, we assume that the Matching Angle Hypothesis is satisfied. 

Consider a sequence $\lambda_l \to \infty$. For each $l$, we consider the domain $\Omega_{\lambda_l}$ defined in Section \ref{smoothing.of.polytope}. Note that $\Omega_{\lambda_l}$ is a compact, convex domain in $\mathbb{R}^n$ with smooth boundary $\partial \Omega_{\lambda_l} = \Sigma_{\lambda_l}$. For each $l$, we define a map $N^{(l)}: \Sigma_{\lambda_l} \to S^{n-1}$ as in Definition \ref{definition.of.N}. Moreover, we define a function $V_{\lambda_l}: \Sigma_{\lambda_l} \to \mathbb{R}$ as in Proposition \ref{formula.for.V}.

Let us fix a Euclidean ball $U$ such that the closure of $U$ is contained in the interior of $\Omega$. Note that $U \subset \Omega_{\lambda_l}$ if $l$ is sufficiently large. In the following, we will always assume that $l$ is chosen sufficiently large so that $U \subset \Omega_{\lambda_l}$.

\begin{proposition}
\label{L2.estimate} 
There exists a uniform constant $C$ (independent of $l$) such that
\[\int_{\Omega_{\lambda_l}} F^2 \, d\sigma_g \leq C \int_{\Omega_{\lambda_l}} |\nabla F|^2 \, d\text{\rm vol}_g + C \int_U F^2 \, d\text{\rm vol}_g\] 
for every smooth function $F: \Omega_{\lambda_l} \to \mathbb{R}$. 
\end{proposition} 

\textbf{Proof.} 
Note that the hypersurface $\Sigma_{\lambda_l} = \partial \Omega_{\lambda_l}$ can be written as a radial graph with bounded slope. From this, it is easy to see that $\Omega_{\lambda_l}$ is bi-Lipschitz equivalent to the Euclidean unit ball, with constants that are independent of $l$. Hence, the assertion follows from the corresponding estimate on the unit ball (which, in turn, is a consequence of the Poincar\'e inequality on the unit ball). \\

\begin{proposition}
\label{sobolev.trace.estimate} 
There exists a uniform constant $C$ (independent of $l$) such that
\[\int_{\Sigma_{\lambda_l}} F^2 \, d\sigma_g \leq C \int_{\Omega_{\lambda_l}} |\nabla F|^2 \, d\text{\rm vol}_g + C \int_{\Omega_{\lambda_l}} F^2 \, d\text{\rm vol}_g\] 
for every smooth function $F: \Omega_{\lambda_l} \to \mathbb{R}$. 
\end{proposition} 

\textbf{Proof.} 
Note that the hypersurface $\Sigma_{\lambda_l} = \partial \Omega_{\lambda_l}$ can be written as a radial graph with bounded slope. From this, it is easy to see that $\Omega_{\lambda_l}$ is bi-Lipschitz equivalent to the Euclidean unit ball, with constants that are independent of $l$. The assertion follows now from the Sobolev trace theorem on the unit ball. \\

\begin{proposition}
\label{fefferman.phong}
We have 
\[\int_{\Sigma_{\lambda_l}} \max \{-V_{\lambda_l},0\} \, F^2 \, d\sigma_g \leq o(1) \int_{\Omega_{\lambda_l}} |\nabla F|^2 \, d\text{\rm vol}_g + o(1) \int_{\Sigma_{\lambda_l}} F^2 \, d\sigma_g\] 
for every smooth function $F: \Omega_{\lambda_l} \to \mathbb{R}$. 
\end{proposition}

\textbf{Proof.} 
Note that the hypersurface $\Sigma_{\lambda_l} = \partial \Omega_{\lambda_l}$ can be written as a radial graph with bounded slope. From this, it is easy to see that $\Omega_{\lambda_l}$ is bi-Lipschitz equivalent to the Euclidean unit ball, with constants that are independent of $l$. Hence, the assertion follows from the Fefferman-Phong estimate on the unit ball (see Corollary \ref{fefferman.phong.ball}) together with Corollary \ref{negative.part.of.V}. This completes the proof of Proposition \ref{fefferman.phong}. \\

For each $l$, we may use the map $N^{(l)}: \Sigma_{\lambda_l} \to S^{n-1}$ to define a boundary chirality $\chi^{(l)}$ (see Definition \ref{definition.of.chi}). It follows from Lemma \ref{homotopy.class.of.N} that $N^{(l)}$ is homotopic to the Euclidean Gauss map of $\Sigma_{\lambda_l}$. By Proposition \ref{index}, we can find an $m$-tuple of spinors $s^{(l)} = (s_1^{(l)},\hdots,s_m^{(l)})$ defined on $\Omega_{\lambda_l}$ with the following properties: 
\begin{itemize} 
\item $s^{(l)}$ is harmonic, i.e. $\mathcal{D} s^{(l)} = 0$ at each point in $\Omega_{\lambda_l}$.
\item $\chi^{(l)} s^{(l)} = s^{(l)}$ at each point on $\Sigma_{\lambda_l}$. 
\item $s^{(l)}$ does not vanish identically. 
\end{itemize} 
Standard unique continuation arguments imply that $\int_U \sum_{\alpha=1}^m |s_\alpha^{(l)}|^2 \, d\text{\rm vol}_g > 0$ if $l$ is sufficiently large. By scaling, we can arrange that $\int_U \sum_{\alpha=1}^m |s_\alpha^{(l)}|^2 \, d\text{\rm vol}_g = m \, \text{\rm vol}_g(U)$ if $l$ is sufficiently large.

\begin{proposition}
\label{covariant.derivative.of.s.small.in.L2}
We have 
\[\int_{\Omega_{\lambda_l}} \sum_{\alpha=1}^m |\nabla s_\alpha^{(l)}|^2 \, d\text{\rm vol}_g \to 0\] 
as $l \to \infty$.
\end{proposition}

\textbf{Proof.} 
Combining Proposition \ref{L2.estimate}, Proposition \ref{sobolev.trace.estimate}, and Proposition \ref{fefferman.phong}, we obtain 
\[\int_{\Sigma_{\lambda_l}} \max \{-V_{\lambda_l},0\} \, F^2 \, d\sigma_g \leq o(1) \int_{\Omega_{\lambda_l}} |\nabla F|^2 \, d\text{\rm vol}_g + o(1) \int_U F^2 \, d\text{\rm vol}_g\] 
for every smooth function $F: \Omega_{\lambda_l} \to \mathbb{R}$. In the next step, we put $F = \big ( \delta^2 + \sum_{\alpha=1}^m |s_\alpha^{(l)}|^2 \big )^{\frac{1}{2}}$, and send $\delta \to 0$. This gives 
\begin{align*} 
&\int_{\Sigma_{\lambda_l}} \max \{-V_{\lambda_l},0\} \, \bigg ( \sum_{\alpha=1}^m \, |s_\alpha^{(l)}|^2 \bigg ) \, d\sigma_g \\ 
&\leq o(1) \int_{\Omega_{\lambda_l}} \sum_{\alpha=1}^m |\nabla s_\alpha^{(l)}|^2 \, d\text{\rm vol}_g + o(1) \int_U \sum_{\alpha=1}^m |s_\alpha^{(l)}|^2 \, d\text{\rm vol}_g. 
\end{align*} 
On the other hand, Proposition \ref{main.inequality} implies 
\begin{align*} 
&\int_{\Omega_{\lambda_l}} \sum_{\alpha=1}^m |\nabla s_\alpha^{(l)}|^2 \, d\text{\rm vol}_g + \frac{1}{4} \int_{\Omega_{\lambda_l}} \sum_{\alpha=1}^m R \, |s_\alpha^{(l)}|^2 \, d\text{\rm vol}_g \\ 
&\leq -\frac{1}{2} \int_{\Sigma_{\lambda_l}} (H-\|dN^{(l)}\|_{\text{\rm tr}}) \, \bigg ( \sum_{\alpha=1}^m \, |s_\alpha^{(l)}|^2 \bigg ) \, d\sigma_g \\ 
&\leq \frac{1}{2} \int_{\Sigma_{\lambda_l}} \max \{-V_{\lambda_l},0\} \, \bigg ( \sum_{\alpha=1}^m \, |s_\alpha^{(l)}|^2 \bigg ) \, d\sigma_g. 
\end{align*} 
Putting these facts together, we conclude that
\begin{align*} 
&\int_{\Omega_{\lambda_l}} \sum_{\alpha=1}^m |\nabla s_\alpha^{(l)}|^2 \, d\text{\rm vol}_g + \frac{1}{4} \int_{\Omega_{\lambda_l}} \sum_{\alpha=1}^m R \, |s_\alpha^{(l)}|^2 \, d\text{\rm vol}_g \\ 
&\leq o(1) \int_{\Omega_{\lambda_l}} \sum_{\alpha=1}^m |\nabla s_\alpha^{(l)}|^2 \, d\text{\rm vol}_g + o(1) \int_U \sum_{\alpha=1}^m |s_\alpha^{(l)}|^2 \, d\text{\rm vol}_g. 
\end{align*} 
By assumption, the scalar curvature of $g$ is nonnegative. If $l$ is sufficiently large, then the first term on the right hand side can be absorbed into the left hand side. This completes the proof of Proposition \ref{covariant.derivative.of.s.small.in.L2}. \\

\begin{corollary} 
\label{L2.bound}
We have 
\[\int_{\Omega_{\lambda_l}} \sum_{\alpha=1}^m |s_\alpha^{(l)}|^2 \, d\text{\rm vol}_g \leq C,\]
where $C$ is a constant that does not depend on $l$. 
\end{corollary} 

\textbf{Proof.} 
We apply Proposition \ref{L2.estimate} with $F = \big ( \delta^2 + \sum_{\alpha=1}^m |s_\alpha^{(l)}|^2 \big )^{\frac{1}{2}}$, and send $\delta \to 0$. This gives 
\[\int_{\Omega_{\lambda_l}} \sum_{\alpha=1}^m |s_\alpha^{(l)}|^2 \, d\text{\rm vol}_g \leq C \int_{\Omega_{\lambda_l}} \sum_{\alpha=1}^m |\nabla s_\alpha^{(l)}|^2 \, d\text{\rm vol}_g + C \int_U \sum_{\alpha=1}^m |s_\alpha^{(l)}|^2 \, d\text{\rm vol}_g,\] 
where $C$ is independent of $l$. Hence, the assertion follows from Proposition \ref{covariant.derivative.of.s.small.in.L2}. This completes the proof of Corollary \ref{L2.bound}. \\

Combining Corollary \ref{L2.bound} with standard interior estimates for elliptic PDE, we obtain smooth estimates for $s^{(l)}$ on compact subsets of $\Omega \setminus \partial \Omega$. After passing to a subsequence if necessary, the sequence $s^{(l)} = (s_1^{(l)},\hdots,s_m^{(l)})$ converges in $C_{\text{\rm loc}}^\infty(\Omega \setminus \partial \Omega)$ to an $m$-tuple of parallel spinors $s = (s_1,\hdots,s_m)$ which is defined on $\Omega \setminus \partial \Omega$. 

\begin{lemma}
\label{boundary.behavior.of.s.1}
We have 
\[\int_{\Sigma_{\lambda_l}} \sum_{\alpha=1}^m |s_\alpha^{(l)}-s_\alpha|^2 \to 0\] 
as $l \to \infty$. 
\end{lemma}

\textbf{Proof.} 
Combining Proposition \ref{L2.estimate} and Proposition \ref{sobolev.trace.estimate}, we obtain 
\[\int_{\Sigma_{\lambda_l}} F^2 \, d\sigma_g \leq C \int_{\Omega_{\lambda_l}} |\nabla F|^2 \, d\text{\rm vol}_g + C \int_U F^2 \, d\text{\rm vol}_g\] 
for every smooth function $F: \Omega_{\lambda_l} \to \mathbb{R}$. In the next step, we put $F = \big ( \delta^2 + \sum_{\alpha=1}^m |s_\alpha^{(l)} - s_\alpha|^2 \big )^{\frac{1}{2}}$, and send $\delta \to 0$. This gives 
\begin{align*} 
&\int_{\Sigma_{\lambda_l}} \sum_{\alpha=1}^m |s_\alpha^{(l)} - s_\alpha|^2 \, d\sigma_g \\ 
&\leq C \int_{\Omega_{\lambda_l}} \sum_{\alpha=1}^m |\nabla (s_\alpha^{(l)} - s_\alpha)|^2 \, d\text{\rm vol}_g + C \int_U \sum_{\alpha=1}^m |s_\alpha^{(l)} - s_\alpha|^2 \, d\text{\rm vol}_g, 
\end{align*} 
where $C$ is independent of $l$. Recall that $s_\alpha$ is parallel for each $\alpha=1,\hdots,m$. Hence, Proposition \ref{covariant.derivative.of.s.small.in.L2} implies that  
\[\int_{\Omega_{\lambda_l}} \sum_{\alpha=1}^m |\nabla (s_\alpha^{(l)} - s_\alpha)|^2 \, d\text{\rm vol}_g = \int_{\Omega_{\lambda_l}} \sum_{\alpha=1}^m |\nabla s_\alpha^{(l)}|^2 \, d\text{\rm vol}_g \to 0\] 
as $l \to \infty$. Moreover, since $s_\alpha^{(l)} \to s_\alpha$ in $C_{\text{\rm loc}}^\infty(\Omega \setminus \partial \Omega)$ for each $\alpha=1,\hdots,m$, we know that 
\[\int_U \sum_{\alpha=1}^m |s_\alpha^{(l)} - s_\alpha|^2 \, d\text{\rm vol}_g \to 0\] 
as $l \to \infty$. This completes the proof of Lemma \ref{boundary.behavior.of.s.1}. \\

\begin{lemma}
\label{boundary.behavior.of.s.2}
We have 
\[\int_{\Sigma_{\lambda_l}} |s-\chi^{(l)} s|^2 \to 0\] 
as $l \to \infty$. 
\end{lemma} 

\textbf{Proof.} 
Recall that $s^{(l)}$ satisfies the boundary condition $\chi^{(l)} s^{(l)} = s^{(l)}$ at each point on $\Sigma_{\lambda_l}$. This implies 
\[|s-\chi^{(l)} s| = |(s-s^{(l)}) - \chi^{(l)} (s-s^{(l)})| \leq C \, |s-s^{(l)}|\] 
at each point on $\Sigma_{\lambda_l}$. Hence, the assertion follows from Lemma \ref{boundary.behavior.of.s.1}. \\

We next analyze the behavior of the map $N^{(l)}: \Sigma_{\lambda_l} \to S^{n-1}$ near the boundary faces of $\Omega$.

\begin{lemma}
\label{boundary.behavior.of.N}
Let us fix an arbitrary element $i_0 \in I$. Suppose that $p \in \{u_{i_0} = 0\} \cap \bigcap_{i \in I \setminus \{i_0\}} \{u_i < 0\}$. Then we can find a small positive real number $r$ (depending on $p$) such that 
\[\sup_{\Sigma_{\lambda_l} \cap B_r(p)} |N^{(l)}-N_{i_0}| \to 0\] 
as $l \to \infty$.
\end{lemma}

\textbf{Proof.} 
This follows directly from the definition of the map $N^{(l)}: \Sigma_{\lambda_l} \to S^{n-1}$. \\

We now continue with the proof of Theorem \ref{main.theorem}. Since $s_1,\hdots,s_m$ are parallel spinors, we can find a fixed matrix $z \in \text{\rm End}(\mathbb{C}^m)$ such that $z_{\alpha\beta} = \langle s_\alpha,s_\beta \rangle$ at each point in the interior of $\Omega$.

\begin{lemma}
\label{z}
For each $i \in I$, the matrix $z \in \text{\rm End}(\mathbb{C}^m)$ commutes with the matrix $\sum_{a=1}^n \langle N_i,E_a \rangle \, \omega_a \in \text{\rm End}(\mathbb{C}^m)$. 
\end{lemma} 

\textbf{Proof.} 
At each point on $\Sigma_{\lambda_l}$, we have 
\begin{align*} 
&\sum_{a=1}^n \sum_{\beta=1}^m \langle N^{(l)},E_a \rangle \, \omega_{a\alpha\beta} \, z_{\beta\gamma} - \sum_{a=1}^n \sum_{\beta=1}^m \langle N^{(l)},E_a \rangle \, \omega_{a\beta\gamma} \, z_{\alpha\beta} \\ 
&= \sum_{a=1}^n \sum_{\beta=1}^m \langle N^{(l)},E_a \rangle \, \omega_{a\alpha\beta} \, \langle \nu \cdot s_\beta,\nu \cdot s_\gamma \rangle + \sum_{a=1}^n \sum_{\beta=1}^m \langle N^{(l)},E_a \rangle \, \overline{\omega_{a\gamma\beta}} \, \langle \nu \cdot s_\alpha,\nu \cdot s_\beta \rangle \\ 
&= -\langle (\chi^{(l)} s)_\alpha,\nu \cdot s_\gamma \rangle - \langle \nu \cdot s_\alpha,(\chi^{(l)} s)_\gamma \rangle \\ 
&= \langle s_\alpha - (\chi^{(l)} s)_\alpha,\nu \cdot s_\gamma \rangle + \langle \nu \cdot s_\alpha,s_\gamma - (\chi^{(l)} s)_\gamma \rangle
\end{align*} 
for all $\alpha,\gamma=1,\hdots,m$. This implies 
\[\bigg | \bigg ( \sum_{a=1}^n \langle N^{(l)},E_a \rangle \, \omega_a \bigg ) \, z - z \, \bigg ( \sum_{a=1}^n \langle N^{(l)},E_a \rangle \, \omega_a \bigg ) \bigg | \leq C \, |s| \, |s - \chi^{(l)} s|\] 
at each point on $\Sigma_{\lambda_l}$. Using Lemma \ref{boundary.behavior.of.s.2}, we conclude that 
\[\int_{\Sigma_{\lambda_l}} \bigg | \bigg ( \sum_{a=1}^n \langle N^{(l)},E_a \rangle \, \omega_a \bigg ) \, z - z \, \bigg ( \sum_{a=1}^n \langle N^{(l)},E_a \rangle \, \omega_a \bigg ) \bigg | \to 0\] 
as $l \to \infty$. We now fix an arbitrary element $i_0 \in I$. By Lemma \ref{boundary.faces}, the set $\{u_{i_0} = 0\} \cap \bigcap_{i \in I \setminus \{i_0\}} \{u_i < 0\}$ is non-empty. Using Lemma \ref{boundary.behavior.of.N}, we deduce that 
\[\bigg ( \sum_{a=1}^n \langle N_{i_0},E_a \rangle \, \omega_a \bigg ) \, z - z \, \bigg ( \sum_{a=1}^n \langle N_{i_0},E_a \rangle \, \omega_a \bigg ) = 0.\] 
This completes the proof of Lemma \ref{z}. \\

Combining Lemma \ref{z} and Lemma \ref{span}, we conclude that the matrix $z \in \text{\rm End}(\mathbb{C}^m)$ commutes with the matrix $\omega_a \in \text{\rm End}(\mathbb{C}^m)$ for each $a=1,\hdots,n$. In view of (\ref{surjectivity.of.hat.rho}), it follows that the matrix $z$ commutes with every element of $\text{\rm End}(\mathbb{C}^m)$. This implies that $z$ is a scalar multiple of the identity. Using the normalization $\int_U \sum_{\alpha=1}^m |s_\alpha|^2 \, d\text{\rm vol}_g = m \, \text{\rm vol}_g(U)$, we conclude that $z$ is the identity.

To summarize, $s = (s_1,\hdots,s_m)$ is a collection of parallel spinors which are defined at each point in the interior of $\Omega$ and are orthonormal at each point in the interior of $\Omega$. Therefore, the spinor bundle is a flat bundle. This implies that the Riemann curvature tensor of $g$ vanishes identically. 


It remains to show that the boundary faces of $\Omega$ are totally geodesic. Recall that $s = (s_1,\hdots,s_m)$ is defined in the interior of $\Omega$. Since $s$ is parallel, we may extend $s$ continuously to $\Omega$. Let us fix an arbitrary element $i_0 \in I$. By Lemma \ref{boundary.faces}, the set $\{u_{i_0} = 0\} \cap \bigcap_{i \in I \setminus \{i_0\}} \{u_i < 0\}$ is non-empty. Using Lemma \ref{boundary.behavior.of.s.2} and Lemma \ref{boundary.behavior.of.N}, we conclude that 
\[\nu_{i_0} \cdot s_\alpha = \sum_{a=1}^n \sum_{\beta=1}^m \langle N_{i_0},E_a \rangle \, \omega_{a\alpha\beta} \, s_\beta\] 
at each point in $\{u_{i_0} = 0\} \cap \bigcap_{i \in I \setminus \{i_0\}} \{u_i < 0\}$. Hence, if $X$ is an arbitrary vector field on $\Omega$, then we obtain 
\begin{align*} 
m \, \langle X,\nu_{i_0} \rangle 
&= \sum_{\alpha=1}^m \langle X,\nu_{i_0} \rangle \, \langle s_\alpha,s_\alpha \rangle \\ 
&= -\frac{1}{2} \sum_{\alpha=1}^m \langle X \cdot \nu_{i_0} \cdot s_\alpha,s_\alpha \rangle + \frac{1}{2} \sum_{\alpha=1}^m \langle X \cdot s_\alpha,\nu_{i_0} \cdot s_\alpha \rangle \\
&= -\frac{1}{2} \sum_{a=1}^n \sum_{\alpha,\beta=1}^m \langle N_{i_0},E_a \rangle \, \omega_{a\alpha\beta} \, \langle X \cdot s_\beta,s_\alpha \rangle \\ 
&+ \frac{1}{2} \sum_{a=1}^n \sum_{\alpha,\beta=1}^m \langle N_{i_0},E_a \rangle \, \overline{\omega_{a\alpha\beta}} \, \langle X \cdot s_\alpha,s_\beta \rangle \\ 
&= -\sum_{a=1}^n \sum_{\alpha,\beta=1}^m \langle N_{i_0},E_a \rangle \, \omega_{a\alpha\beta} \, \langle X \cdot s_\beta,s_\alpha \rangle
\end{align*} 
at each point in $\{u_{i_0} = 0\} \cap \bigcap_{i \in I \setminus \{i_0\}} \{u_i < 0\}$. Since $s = (s_1,\hdots,s_m)$ is parallel, it follows that $\nu_{i_0}$ is parallel along the hypersurface $\{u_{i_0} = 0\} \cap \bigcap_{i \in I \setminus \{i_0\}} \{u_i < 0\}$. Consequently, the second fundamental form of the hypersurface $\{u_{i_0} = 0\}$ vanishes at each point in $\{u_{i_0} = 0\} \cap \bigcap_{i \in I \setminus \{i_0\}} \{u_i < 0\}$. In view of Lemma \ref{boundary.faces}, we conclude that the second fundamental form of the hypersurface $\{u_{i_0} = 0\}$ vanishes at each point in $\Omega \cap \{u_{i_0} = 0\}$.

\appendix

\section{A variant of a theorem of Fefferman and Phong} 

In this section, we describe a variant of an estimate due to Fefferman and Phong \cite{Fefferman-Phong}, which plays a central role in our argument. Throughout this section, we fix an integer $n \geq 3$. We denote by $\mathcal{Q}$ the collection of all $(n-1)$-dimensional cubes of the form 
\[[2^m j_1,2^m (j_1+1)] \times \hdots \times [2^m j_{n-1},2^m (j_{n-1}+1)] \times \{0\},\] 
where $m \in \mathbb{Z}$ and $j_1,\hdots,j_{n-1} \in \mathbb{Z}$. For abbreviation, we put 
\[\Gamma = \bigcup_{Q \in \mathcal{Q}} \partial Q.\] 
If $Q_1,Q_2 \in \mathcal{Q}$ satisfy $Q_1 \cap Q_2 \setminus \Gamma \neq \emptyset$, then $Q_1 \subset Q_2$ or $Q_2 \subset Q_1$.

For each $(n-1)$-dimensional cube $Q \in \mathcal{Q}$, we denote by $|Q|$ the $(n-1)$-dimensional volume of $Q$. 

\begin{theorem}
\label{fefferman.phong.halfspace}
Let us fix an integer $n \geq 3$ and a real number $\sigma \in (1,n-1)$. Suppose that $V$ is a nonnegative continuous function defined on the hyperplane $\mathbb{R}^{n-1} \times \{0\}$ with the property that 
\begin{equation} 
\label{assumption.V.halfspace}
\bigg ( |Q|^{-1} \int_Q V^\sigma \bigg )^{\frac{1}{\sigma}} \leq \text{\rm diam}(Q)^{-1} 
\end{equation}
for each $(n-1)$-dimensional cube $Q \in \mathcal{Q}$. Let $F$ be a smooth function defined on the half-space $\mathbb{R}_+^n = \{x \in \mathbb{R}^n: x_n \geq 0\}$, and let $f$ denote the restriction of $F$ to the boundary $\partial \mathbb{R}_+^n = \mathbb{R}^{n-1} \times \{0\}$. Then 
\[\int_Q V f^2 \leq C \int_{Q \times [0,\text{\rm diam}(Q)]} |\nabla F|^2 + C \, \text{\rm diam}(Q)^{-1} \int_Q f^2.\] 
for each $(n-1)$-dimensional cube $Q \in \mathcal{Q}$. The constant $C$ depends only on $n$ and $\sigma$.
\end{theorem}

The proof of Theorem \ref{fefferman.phong.halfspace} is a straightforward adaptation of the arguments of Fefferman and Phong \cite{Fefferman-Phong}. Let us fix an exponent $\tau \in (1,\sigma)$. Let $V: \mathbb{R}^{n-1} \times \{0\} \to \mathbb{R}$ be a nonnegative continuous function satisfying (\ref{assumption.V.halfspace}). We define a measurable function $W: \mathbb{R}^{n-1} \times \{0\} \to \mathbb{R}$ by 
\[W(x) = \sup_{Q \in \mathcal{Q}, x \in Q} \bigg ( |Q|^{-1} \int_Q V^\sigma \bigg )^{\frac{1}{\sigma}}\] 
for each point $x \in \mathbb{R}^{n-1} \times \{0\}$. It follows from (\ref{assumption.V.halfspace}) that $W$ is locally bounded. Moreover, $V \leq W$ at each point in $\mathbb{R}^{n-1} \times \{0\}$. 

Let $F$ be a smooth function defined on the half-space $\mathbb{R}_+^n = \{x \in \mathbb{R}^n: x_n \geq 0\}$, and let $f$ denote the restriction of $F$ to the boundary $\partial \mathbb{R}_+^n = \mathbb{R}^{n-1} \times \{0\}$. For each $(n-1)$-dimensional cube $Q \in \mathcal{Q}$, we denote by $f_Q = |Q|^{-1} \int_Q f$ the mean value of $f$ over the cube $Q$. 

\begin{lemma}
\label{estimate.for.W.in.L.tau}
For each $(n-1)$-dimensional cube $Q_0 \in \mathcal{Q}$, we have 
\[\bigg ( |Q_0|^{-1} \int_{Q_0} W^\tau \bigg )^{\frac{1}{\tau}} \leq C \, \sup_{Q \in \mathcal{Q}, Q_0 \subset Q} \bigg ( |Q|^{-1} \int_Q V^\sigma \bigg )^{\frac{1}{\sigma}},\] 
where $C$ depends only on $n$, $\sigma$, and $\tau$.
\end{lemma}

\textbf{Proof.} 
For abbreviation, let 
\[\Lambda = \sup_{Q \in \mathcal{Q}, Q_0 \subset Q} \bigg ( |Q|^{-1} \int_Q V^\sigma \bigg )^{\frac{1}{\sigma}}.\] 
It follows from (\ref{assumption.V.halfspace}) that $\Lambda<\infty$. We define a bounded measurable function $W_0: Q_0 \to \mathbb{R}$ by 
\[W_0(x) = \sup_{Q \in \mathcal{Q}, x \in Q \subset Q_0} \bigg ( |Q|^{-1} \int_Q V^\sigma \bigg )^{\frac{1}{\sigma}}\] 
for each point $x \in Q_0$. Then 
\[W(x) = \max \{\Lambda,W_0(x)\}\] 
for each point $x \in Q_0 \setminus \Gamma$. The function $W_0^\sigma$ is bounded from above by the maximal function associated with the function $V^\sigma \, 1_{Q_0}$. Hence, the weak version of the Hardy-Littlewood maximal inequality (cf. \cite{Muscalu-Schlag}, Proposition 2.9 (i)) implies 
\begin{equation} 
\label{inequality.for.level.set.W_0}
|Q_0|^{-1} \, |\{x \in Q_0: W_0(x)^\sigma > \alpha\}| \leq C \alpha^{-1} \, |Q_0|^{-1} \int_{Q_0} V^\sigma \leq C \alpha^{-1} \, \Lambda^\sigma 
\end{equation}
for all $\alpha>0$. We multiply both sides of (\ref{inequality.for.level.set.W_0}) by $\frac{\tau}{\sigma} \, \alpha^{\frac{\tau}{\sigma}-1}$ and integrate over $\alpha \in (\Lambda^\sigma,\infty)$. Using Fubini's theorem, we obtain 
\begin{align*}
&|Q_0|^{-1} \int_{Q_0} \max \{W_0^\tau-\Lambda^\tau,0\} \\ 
&= |Q_0|^{-1} \int_{\Lambda^\sigma}^\infty \frac{\tau}{\sigma} \, \alpha^{\frac{\tau}{\sigma}-1} \, |\{x \in Q_0: W_0(x)^\sigma > \alpha\}| \, d\alpha \\ 
&\leq C \int_{\Lambda^\sigma}^\infty \frac{\tau}{\sigma} \, \alpha^{\frac{\tau}{\sigma}-2} \, \Lambda^\sigma \, d\alpha \\ 
&= C \, \frac{\tau}{\sigma-\tau} \, \Lambda^\tau. 
\end{align*} 
Putting these facts together, we conclude that 
\[|Q_0|^{-1} \int_{Q_0} W^\tau \leq C \, \Lambda^\tau.\] 
This completes the proof of Lemma \ref{estimate.for.W.in.L.tau}. \\

\begin{lemma}
\label{consequence.of.reverse.Holder.inequality}
Given a real number $\varepsilon>0$, we can find a real number $\delta>0$ (depending only on $n$, $\sigma$, and $\varepsilon$) with the property that 
\[\int_A W \leq \varepsilon \int_{Q_0} W.\] 
for every $(n-1)$-dimensional cube $Q_0 \in \mathcal{Q}$ and every measurable set $A \subset Q_0$ satisfying $|A| \leq \delta \, |Q_0|$.
\end{lemma}

\textbf{Proof.} 
Using Lemma \ref{estimate.for.W.in.L.tau}, we obtain 
\[\bigg ( |Q_0|^{-1} \int_{Q_0} W^\tau \bigg )^{\frac{1}{\tau}} \leq C \, \sup_{Q \in \mathcal{Q}, Q_0 \subset Q} \bigg ( |Q|^{-1} \int_Q V^\sigma \bigg )^{\frac{1}{\sigma}}.\] 
Moreover, 
\[\sup_{Q \in \mathcal{Q}, Q_0 \subset Q} \bigg ( |Q|^{-1} \int_Q V^\sigma \bigg )^{\frac{1}{\sigma}} \leq \inf_{Q_0} W \leq |Q_0|^{-1} \int_{Q_0} W\] 
by definition of $W$. Putting these facts together, we obtain 
\[\bigg ( |Q_0|^{-1} \int_{Q_0} W^\tau \bigg )^{\frac{1}{\tau}} \leq C \, |Q_0|^{-1} \int_{Q_0} W.\] 
Hence, if $A \subset Q_0$ is a measurable set with $|A| \leq \delta \, |Q_0|$, then H\"older's inequality gives 
\[\int_A W \leq |A|^{\frac{\tau-1}{\tau}} \, \bigg ( \int_{Q_0} W^\tau \bigg )^{\frac{1}{\tau}} \leq \delta^{\frac{\tau-1}{\tau}} \, |Q_0|^{\frac{\tau-1}{\tau}} \, \bigg ( \int_{Q_0} W^\tau \bigg )^{\frac{1}{\tau}} \leq C \delta^{\frac{\tau-1}{\tau}} \int_{Q_0} W.\] 
Hence, if we choose $\delta$ to be a small multiple of $\varepsilon^{\frac{\tau}{\tau-1}}$, then $\delta$ has the required property. This completes the proof of Lemma \ref{consequence.of.reverse.Holder.inequality}. \\

\begin{lemma} 
\label{estimate.for.W} 
For each $(n-1)$-dimensional cube $Q_0 \in \mathcal{Q}$, we have 
\[|Q_0|^{-1} \int_{Q_0} W \leq C \, \text{\rm diam}(Q_0)^{-1},\] 
where $C$ depends only on $n$ and $\sigma$.
\end{lemma} 

\textbf{Proof.} 
Using Lemma \ref{estimate.for.W.in.L.tau} and H\"older's inequality, we obtain 
\[|Q_0|^{-1} \int_{Q_0} W \leq \bigg ( |Q_0|^{-1} \int_{Q_0} W^\tau \bigg )^{\frac{1}{\tau}} \leq C \, \sup_{Q \in \mathcal{Q}, Q_0 \subset Q} \bigg ( |Q|^{-1} \int_Q V^\sigma \bigg )^{\frac{1}{\sigma}}.\] 
Hence, the assertion follows from (\ref{assumption.V.halfspace}). \\

\begin{lemma}
\label{main.estimate.1}
Let us fix an $(n-1)$-dimensional cube $Q_0 \in \mathcal{Q}$. We define a bounded measurable function $g: Q_0 \to \mathbb{R}$ by 
\[g(x) = \sup_{Q \in \mathcal{Q}, x \in Q \subset Q_0} |Q|^{-1} \int_Q |f-f_Q|\] 
for each point $x \in Q_0$. Then 
\[\int_{Q_0} V \, |f-f_{Q_0}|^2 \leq C \int_{Q_0} W g^2,\] 
where $C$ depends only on $n$ and $\sigma$.
\end{lemma}

\textbf{Proof.} 
We define a bounded measurable function $h: Q_0 \to \mathbb{R}$ by 
\[h(x) = \sup_{Q \in \mathcal{Q}, x \in Q \subset Q_0} |Q|^{-1} \int_Q |f-f_{Q_0}|\] 
for each point $x \in Q_0$. Note that $V \leq W$ and $|f-f_{Q_0}| \leq h$ at each point in $Q_0$. Hence, it suffices to prove that 
\begin{equation} 
\label{important.estimate}
\int_{Q_0} W h^2 \leq C \int_{Q_0} W g^2.
\end{equation} 
In order to prove the inequality (\ref{important.estimate}), we define $\alpha_0 = |Q_0|^{-1} \int_{Q_0} |f-f_{Q_0}|$. For each $\alpha > \alpha_0$, we denote by $\mathcal{Q}_\alpha$ the set of all $(n-1)$-dimensional cubes $Q \in \mathcal{Q}$ with the following properties: 
\begin{itemize}
\item $Q \subset Q_0$. 
\item $|Q|^{-1} \int_Q |f-f_{Q_0}| > \alpha$.
\item If $\tilde{Q} \in \mathcal{Q}$ is an $(n-1)$-dimensional cube with $Q \subsetneq \tilde{Q}$ and $\tilde{Q} \subset Q_0$, then $|\tilde{Q}|^{-1} \int_{\tilde{Q}} |f-f_{Q_0}| \leq \alpha$.
\end{itemize} 
It follows from the definition of $\alpha_0$ that $Q_0 \notin \mathcal{Q}_\alpha$ for each $\alpha>\alpha_0$. It is easy to see that 
\begin{equation} 
\label{property.of.Q.1}
\alpha < |Q|^{-1} \int_Q |f-f_{Q_0}| \leq 2^{n-1} \alpha 
\end{equation}
for each $\alpha>\alpha_0$ and each $Q \in \mathcal{Q}_\alpha$. Moreover, 
\begin{equation} 
\{h>\alpha\} = \bigcup_{Q \in \mathcal{Q}_\alpha} Q 
\end{equation}
for each $\alpha>\alpha_0$. Finally, given a real number $\alpha>\alpha_0$ and a point $x \notin \Gamma$, there is at most one cube $Q \in \mathcal{Q}_\alpha$ that contains the point $x$.

We next apply Lemma \ref{consequence.of.reverse.Holder.inequality} with $\varepsilon = 2^{-2n-1}$. Hence, we can find a real number $\delta \in (0,1)$ such that 
\begin{equation} 
\label{A}
\int_A W \leq 2^{-2n-1} \int_Q W 
\end{equation}
for every $(n-1)$-dimensional cube $Q \in \mathcal{Q}$ and every measurable set $A \subset Q$ satisfying $|A| \leq 2^{1-n} \, \delta \, |Q|$.

Let us consider a real number $\alpha>\alpha_0$ and an $(n-1)$-dimensional cube $Q \in \mathcal{Q}_\alpha$. The upper bound in (\ref{property.of.Q.1}) implies $|f_Q-f_{Q_0}| \leq 2^{n-1} \alpha$. Using the lower bound in (\ref{property.of.Q.1}), we obtain 
\[2^{n-1} \alpha \, |\tilde{Q}| \leq \int_{\tilde{Q}} (|f-f_{Q_0}| - 2^{n-1} \alpha) \leq \int_{\tilde{Q}} |f-f_Q|\] 
for all $(n-1)$-dimensional cubes $\tilde{Q} \in \mathcal{Q}_{2^n \alpha}$. In the next step, we take the sum over all $(n-1)$-dimensional cubes $\tilde{Q} \in \mathcal{Q}_{2^n \alpha}$ with $\tilde{Q} \subset Q$. This gives 
\begin{equation} 
\label{summation.1}
2^{n-1} \alpha \sum_{\tilde{Q} \in \mathcal{Q}_{2^n \alpha}, \tilde{Q} \subset Q} |\tilde{Q}| \leq \int_Q |f-f_Q|. 
\end{equation}
For each $(n-1)$-dimensional cube $Q \in \mathcal{Q}_\alpha$, we have the inclusion
\begin{equation} 
\label{inclusion.1}
Q \cap \{h > 2^n \alpha\} \setminus \Gamma \subset \bigcup_{\tilde{Q} \in \mathcal{Q}_{2^n \alpha}, \tilde{Q} \subset Q} \tilde{Q}. 
\end{equation}
Combining (\ref{summation.1}) and (\ref{inclusion.1}), we conclude that 
\begin{equation} 
2^{n-1} \alpha \, |Q \cap \{h > 2^n \alpha\}| \leq \int_Q |f-f_Q| 
\end{equation}
for every $(n-1)$-dimensional cube $Q \in \mathcal{Q}_\alpha$. In particular, 
\[|Q \cap \{h > 2^n \alpha\}| \leq 2^{1-n} \, \delta \, |Q|\] 
for every $(n-1)$-dimensional cube $Q \in \mathcal{Q}_\alpha$ satisfying $|Q|^{-1} \int_Q |f-f_Q| \leq \delta \alpha$. Applying (\ref{A}) with $A = Q \cap \{h > 2^n \alpha\}$ gives 
\begin{equation} 
\label{case.1}
\int_{Q \cap \{h > 2^n \alpha\}} W \leq 2^{-2n-1} \int_Q W 
\end{equation} 
for every $(n-1)$-dimensional cube $Q \in \mathcal{Q}_\alpha$ satisfying $|Q|^{-1} \int_Q |f-f_Q| \leq \delta \alpha$. On the other hand, if $Q \in \mathcal{Q}_\alpha$ is an $(n-1)$-dimensional cube satisfying $|Q|^{-1} \int_Q |f-f_Q| > \delta \alpha$, then $g > \delta \alpha$ at each point in $Q$. Therefore, 
\begin{equation} 
\label{case.2}
\int_{Q \cap \{h > 2^n \alpha\}} W \leq \int_{Q \cap \{g > \delta\alpha\}} W 
\end{equation} 
for every $(n-1)$-dimensional cube $Q \in \mathcal{Q}_\alpha$ satisfying $|Q|^{-1} \int_Q |f-f_Q| > \delta \alpha$. Combining (\ref{case.1}) and (\ref{case.2}), we conclude that 
\begin{equation} 
\int_{Q \cap \{h > 2^n \alpha\}} W \leq 2^{-2n-1} \int_Q W + \int_{Q \cap \{g > \delta\alpha\}} W 
\end{equation} 
for every $(n-1)$-dimensional cube $Q \in \mathcal{Q}_\alpha$. Summation over all $(n-1)$-dimensional cubes $Q \in \mathcal{Q}_\alpha$  gives 
\begin{equation} 
\label{inequality.for.level.set.of.h}
\int_{\{h > 2^n \alpha\}} W \leq 2^{-2n-1} \int_{\{h > \alpha\}} W + \int_{\{g > \delta\alpha\}} W 
\end{equation} 
for each $\alpha>\alpha_0$. 

We now multiply both sides of (\ref{inequality.for.level.set.of.h}) by $2\alpha$ and integrate over $\alpha \in (2\alpha_0,\infty)$. Using Fubini's theorem, we obtain 
\begin{align*} 
&\int_{Q_0} W \, \max \{2^{-2n} h^2 - 4\alpha_0^2,0\} \\ 
&= \int_{2\alpha_0}^\infty 2\alpha \, \bigg ( \int_{\{h > 2^n \alpha\}} W \bigg ) \, d\alpha \\ 
&\leq \int_{2\alpha_0}^\infty 2^{-2n} \alpha \, \bigg ( \int_{\{h > \alpha\}} W \bigg ) \, d\alpha + \int_{2\alpha_0}^\infty 2\alpha \, \bigg ( \int_{\{g > \delta\alpha\}} W \bigg ) \, d\alpha \\ 
&\leq 2^{-2n-1} \int_{Q_0} W h^2 + \delta^{-2} \int_{Q_0} W g^2. 
\end{align*} 
Rearranging terms gives  
\begin{equation} 
2^{-2n-1} \int_{Q_0} W h^2 \leq 4\alpha_0^2 \int_{Q_0} W + \delta^{-2} \int_{Q_0} W g^2. 
\end{equation}
On the other hand, $g \geq |Q_0|^{-1} \int_{Q_0} |f-f_{Q_0}| = \alpha_0$ at each point in $Q_0$. Consequently, 
\begin{equation} 
\label{final.estimate}
2^{-2n-1} \int_{Q_0} W h^2 \leq (4+\delta^{-2}) \int_{Q_0} W g^2. 
\end{equation} 
The inequality (\ref{important.estimate}) follows immediately from (\ref{final.estimate}). This completes the proof of Lemma \ref{main.estimate.1}. \\

\begin{lemma} 
\label{main.estimate.2}
Let us fix an $(n-1)$-dimensional cube $Q_0 \in \mathcal{Q}$. We define a bounded measurable function $g: Q_0 \to \mathbb{R}$ by 
\[g(x) = \sup_{Q \in \mathcal{Q}, x \in Q \subset Q_0} |Q|^{-1} \int_Q |f-f_Q|\] 
for each point $x \in Q_0$. Then 
\[\int_{Q_0} W g^2 \leq C \int_{Q_0 \times [0,\text{\rm diam}(Q_0)]} |\nabla F|^2,\] 
where $C$ depends only on $n$ and $\sigma$.
\end{lemma}

\textbf{Proof.} 
Let $\alpha_0 = |Q_0|^{-1} \int_{Q_0} |f-f_{Q_0}|$. For each $\alpha > \alpha_0$, we denote by $\mathcal{Q}_\alpha$ the set of all $(n-1)$-dimensional cubes $Q \in \mathcal{Q}$ with the following properties: 
\begin{itemize}
\item $Q \subset Q_0$. 
\item $|Q|^{-1} \int_Q |f-f_Q| > \alpha$.
\item If $\tilde{Q} \in \mathcal{Q}$ is an $(n-1)$-dimensional cube with $Q \subsetneq \tilde{Q}$ and $\tilde{Q} \subset Q_0$, then $|\tilde{Q}|^{-1} \int_{\tilde{Q}} |f-f_{\tilde{Q}}| \leq \alpha$.
\end{itemize} 
It follows from the definition of $\alpha_0$ that $Q_0 \notin \mathcal{Q}_\alpha$ for each $\alpha>\alpha_0$. It is easy to see that 
\begin{equation} 
\label{property.of.Q.2}
\alpha < |Q|^{-1} \int_Q |f-f_Q| \leq 2^n \alpha 
\end{equation}
for each $\alpha>\alpha_0$ and each $Q \in \mathcal{Q}_\alpha$. Moreover, 
\begin{equation} 
\{g>\alpha\} = \bigcup_{Q \in \mathcal{Q}_\alpha} Q 
\end{equation}
for each $\alpha>\alpha_0$. Finally, given a real number $\alpha>\alpha_0$ and a point $x \notin \Gamma$, there is at most one cube $Q \in \mathcal{Q}_\alpha$ that contains the point $x$.

Let us consider a real number $\alpha>\alpha_0$ and an $(n-1)$-dimensional cube $Q \in \mathcal{Q}_\alpha$. Using the lower bound in (\ref{property.of.Q.2}), we obtain 
\[2^{n+2} \alpha \, |\tilde{Q}| \leq \int_{\tilde{Q}} |f-f_{\tilde{Q}}| \leq 2 \int_{\tilde{Q}} |f-f_Q|\] 
for all $(n-1)$-dimensional cubes $\tilde{Q} \in \mathcal{Q}_{2^{n+2} \alpha}$. In the next step, we take the sum over all $(n-1)$-dimensional cubes $\tilde{Q} \in \mathcal{Q}_{2^{n+2} \alpha}$ with $\tilde{Q} \subset Q$. Using the upper bound in (\ref{property.of.Q.2}), we obtain 
\begin{equation} 
\label{summation.2}
2^{n+2} \alpha \sum_{\tilde{Q} \in \mathcal{Q}_{2^{n+2} \alpha}, \tilde{Q} \subset Q} |\tilde{Q}| \leq 2 \int_Q |f-f_Q| \leq 2^{n+1} \alpha \, |Q| 
\end{equation}
for every $(n-1)$-dimensional cube $Q \in \mathcal{Q}_\alpha$. For each $(n-1)$-dimensional cube $Q \in \mathcal{Q}_\alpha$, we have the inclusion
\begin{equation} 
\label{inclusion.2}
Q \cap \{g > 2^{n+2} \alpha\} \setminus \Gamma \subset \bigcup_{\tilde{Q} \in \mathcal{Q}_{2^{n+2} \alpha}, \tilde{Q} \subset Q} \tilde{Q}. 
\end{equation}
Combining (\ref{summation.2}) and (\ref{inclusion.2}), we conclude that 
\[2^{n+2} \alpha \, |Q \cap \{g>2^{n+2} \alpha\}| \leq 2^{n+1} \alpha \, |Q|,\] 
hence 
\begin{equation} 
\label{volume.of.sub.level.set.of.g}
|Q \cap \{g \leq 2^{n+2} \alpha\}| \geq \frac{1}{2} \, |Q| 
\end{equation}
for every $(n-1)$-dimensional cube $Q \in \mathcal{Q}_\alpha$. 

We define a nonnegative function $\varphi: \mathbb{R}^{n-1} \times \{0\} \to \mathbb{R}$ by 
\begin{equation} 
\label{def.varphi}
\varphi(x_1,\hdots,x_{n-1},0) = \bigg ( \int_0^{\text{\rm diam}(Q_0)} |\nabla F(x_1,\hdots,x_{n-1},x_n)|^2 \, dx_n \bigg )^{\frac{1}{2}}. 
\end{equation}
Moreover, we define a nonnegative function $\psi: Q_0 \to \mathbb{R}$ by 
\begin{equation} 
\label{def.psi}
\psi(x) = \sup_{Q \in \mathcal{Q}, x \in Q \subset Q_0} |Q|^{-1} \int_Q \varphi 
\end{equation}
for each point $x \in Q_0$. Using the Sobolev trace theorem, we obtain 
\begin{align} 
\label{sobolev.trace}
\alpha 
&\leq |Q|^{-1} \int_Q |f-f_Q| \notag \\ 
&\leq 2 \, |Q|^{-1} \inf_{a \in \mathbb{R}} \int_Q |f-a| \notag \\ 
&\leq C \, \, |Q|^{-1} \inf_{a \in \mathbb{R}} \bigg ( \int_{Q \times [0,\text{\rm diam}(Q)]} |\nabla (F-a)| \\ 
&\hspace{30mm} + \text{\rm diam}(Q)^{-1} \int_{Q \times [0,\text{\rm diam}(Q)]} |F-a| \bigg ) \notag
\end{align} 
for every $(n-1)$-dimensional cube $Q \in \mathcal{Q}_\alpha$. Using (\ref{sobolev.trace}) and the Poincar\'e inequality, we conclude that 
\begin{equation} 
\alpha \leq C \, |Q|^{-1} \int_{Q \times [0,\text{\rm diam}(Q)]} |\nabla F| 
\end{equation} 
for every $(n-1)$-dimensional cube $Q \in \mathcal{Q}_\alpha$. Using H\"older's inequality, we deduce that 
\begin{equation} 
\label{lower.bound.for.inf.psi}
\alpha \leq C \, \text{\rm diam}(Q)^{\frac{1}{2}} \, |Q|^{-1} \int_Q \varphi \leq C \, \text{\rm diam}(Q)^{\frac{1}{2}} \, \inf_Q \psi 
\end{equation} 
for every $(n-1)$-dimensional cube $Q \in \mathcal{Q}_\alpha$. Combining (\ref{volume.of.sub.level.set.of.g}) and (\ref{lower.bound.for.inf.psi}) gives 
\begin{align} 
\label{lower.bound.for.integral.of.psi.squared}
\alpha^2 \, \text{\rm diam}(Q)^{-1} \, |Q| 
&\leq 2\alpha^2 \, \text{\rm diam}(Q)^{-1} \, |Q \cap \{g \leq 2^{n+2} \alpha\}| \notag \\
&\leq C \int_{Q \cap \{g \leq 2^{n+2} \alpha\}} \psi^2 
\end{align} 
for every $(n-1)$-dimensional cube $Q \in \mathcal{Q}_\alpha$. Combining the estimate (\ref{lower.bound.for.integral.of.psi.squared}) with Lemma \ref{estimate.for.W}, we obtain 
\begin{equation} 
\alpha^2 \int_Q W \leq C \int_{Q \cap \{g \leq 2^{n+2} \alpha\}} \psi^2 
\end{equation}
for every $(n-1)$-dimensional cube $Q \in \mathcal{Q}_\alpha$. Summation over all $(n-1)$-dimensional cubes $Q \in \mathcal{Q}_\alpha$ gives 
\begin{equation} 
\label{lower.bound.for.integral.of.psi.squared.2}
\alpha^2 \int_{\{g > \alpha\}} W \leq C \int_{\{\alpha < g \leq 2^{n+2} \alpha\}} \psi^2 
\end{equation} 
for each $\alpha>\alpha_0$. 

We now multiply both sides of (\ref{lower.bound.for.integral.of.psi.squared.2}) by $2\alpha^{-1}$ and integrate over $\alpha \in (2\alpha_0,\infty)$. Using Fubini's theorem, we obtain  
\begin{align} 
\label{lower.bound.for.integral.of.psi.squared.3}
\int_{Q_0} W \, \max \{g^2-4\alpha_0^2,0\} 
&= \int_{2\alpha_0}^\infty 2\alpha \, \bigg ( \int_{\{g > \alpha\}} W \bigg ) \, d\alpha \notag \\ 
&\leq \int_{2\alpha_0}^\infty 2C \alpha^{-1} \, \bigg ( \int_{\{\alpha < g \leq 2^{n+2} \alpha\}} \psi^2 \bigg ) \, d\alpha \\ 
&\leq 2C \log (2^{n+2}) \int_{Q_0} \psi^2. \notag
\end{align}
On the other hand, the function $\psi$ is bounded from above by the maximal function associated with the function $\varphi \, 1_{Q_0}$. Hence, the strong version of the Hardy-Littlewood maximal inequality (cf. \cite{Muscalu-Schlag}, Proposition 2.9 (ii)) implies 
\begin{equation} 
\label{application.of.hardy.littlewood}
\int_{Q_0} \psi^2 \leq C \int_{Q_0} \varphi^2 \leq C \int_{Q_0 \times [0,\text{\rm diam}(Q_0)]} |\nabla F|^2. 
\end{equation}
Combining (\ref{lower.bound.for.integral.of.psi.squared.3}) and (\ref{application.of.hardy.littlewood}) gives 
\begin{equation} 
\label{estimate.for.g.large}
\int_{Q_0} W \, \max \{g^2-4\alpha_0^2,0\} \leq C \int_{Q_0 \times [0,\text{\rm diam}(Q_0)]} |\nabla F|^2. 
\end{equation}
Finally, using the Sobolev trace theorem, we obtain 
\begin{align} 
\label{sobolev.trace.2}
\alpha_0 
&= |Q_0|^{-1} \int_{Q_0} |f-f_{Q_0}| \notag \\ 
&\leq 2 \, |Q_0|^{-1} \inf_{a \in \mathbb{R}} \int_{Q_0} |f-a| \notag \\ 
&\leq C \, \, |Q_0|^{-1} \inf_{a \in \mathbb{R}} \bigg ( \int_{Q_0 \times [0,\text{\rm diam}(Q_0)]} |\nabla (F-a)| \\ 
&\hspace{30mm} + \text{\rm diam}(Q_0)^{-1} \int_{Q_0 \times [0,\text{\rm diam}(Q_0)]} |F-a| \bigg ). \notag
\end{align} 
Using (\ref{sobolev.trace.2}) and the Poincar\'e inequality, we conclude that 
\begin{equation} 
\label{bound.for.alpha_0}
\alpha_0 \leq C \, |Q_0|^{-1} \int_{Q_0 \times [0,\text{\rm diam}(Q_0)]} |\nabla F|. 
\end{equation}
Using H\"older's inequality, we deduce that 
\begin{equation}
\label{bound.for.alpha_0.2}
\alpha_0^2 \, \text{\rm diam}(Q_0)^{-1} \, |Q_0| \leq C \int_{Q_0 \times [0,\text{\rm diam}(Q_0)]} |\nabla F|^2. 
\end{equation}
Combing the estimate (\ref{bound.for.alpha_0.2}) with Lemma \ref{estimate.for.W} gives 
\begin{equation}
\label{bound.for.alpha_0.3} 
\alpha_0^2 \int_{Q_0} W \leq C \int_{Q_0 \times [0,\text{\rm diam}(Q_0)]} |\nabla F|^2. 
\end{equation}
The assertion follows by combining (\ref{estimate.for.g.large}) and (\ref{bound.for.alpha_0.3}). This completes the proof of Lemma \ref{main.estimate.2}. \\

After these preparations, we now complete the proof of Theorem \ref{fefferman.phong.halfspace}. Combining Lemma \ref{main.estimate.1} and Lemma \ref{main.estimate.2}, we conclude that 
\[\int_{Q_0} V \, |f-f_{Q_0}|^2 \leq C \int_{Q_0 \times [0,\text{\rm diam}(Q_0)]} |\nabla F|^2\] 
for every $(n-1)$-dimensional cube $Q_0 \in \mathcal{Q}$. This implies 
\[\int_{Q_0} V f^2 \leq C \int_{Q_0 \times [0,\text{\rm diam}(Q_0)]} |\nabla F|^2 + C \, |Q_0|^{-2} \, \bigg ( \int_{Q_0} V \bigg ) \,  \bigg ( \int_{Q_0} |f| \bigg )^2\] 
for each $(n-1)$-dimensional cube $Q_0 \in \mathcal{Q}$. Moreover, 
\[|Q_0|^{-1} \int_{Q_0} V \leq \bigg ( |Q_0|^{-1} \int_{Q_0} V^\sigma \bigg )^{\frac{1}{\sigma}} \leq \text{\rm diam}(Q_0)^{-1}\] 
by (\ref{assumption.V.halfspace}). Thus, we conclude that 
\[\int_{Q_0} V f^2 \leq C \int_{Q_0 \times [0,\text{\rm diam}(Q_0)]} |\nabla F|^2 + C \, \text{\rm diam}(Q_0)^{-1} \, |Q_0|^{-1} \, \bigg ( \int_{Q_0} |f| \bigg )^2\] 
for each $(n-1)$-dimensional cube $Q_0 \in \mathcal{Q}$. This completes the proof of Theorem \ref{fefferman.phong.halfspace}. \\

\begin{corollary}
\label{fefferman.phong.ball}
Let us fix an integer $n \geq 3$ and a real number $\sigma \in (1,n-1)$. Suppose that $V$ is a nonnegative continuous function defined on the unit sphere $S^{n-1} \subset \mathbb{R}^n$ with the property that 
\begin{equation} 
\label{assumption.V.ball}
\bigg ( r^{\sigma+1-n} \int_{S^{n-1} \cap B_r(p)} V^\sigma \bigg )^{\frac{1}{\sigma}} \leq 1 
\end{equation}
for all points $p \in \mathbb{R}^n$ and all $0 < r \leq 1$. Let $F$ be a smooth function defined on the unit ball $B^n = \{x \in \mathbb{R}^n: |x| \leq 1\}$, and let $f$ denote the restriction of $F$ to the boundary $\partial B^n = S^{n-1}$. Then 
\[\int_{S^{n-1}} V f^2 \leq C \int_{B^n} |\nabla F|^2 + C \int_{S^{n-1}} f^2.\] 
The constant $C$ depends only on $n$ and $\sigma$.
\end{corollary}

\end{document}